\documentclass[11pt, a4paper, reqno]{amsart}

\usepackage{amssymb}
\usepackage{mathrsfs}
\usepackage{textcomp}
\usepackage{graphicx}

\usepackage[breaklinks]{hyperref}

\hypersetup{
unicode=true,
colorlinks=true,
linkcolor=blue,
citecolor=blue,
urlcolor=blue,
filecolor=blue,
bookmarksnumbered=true,
pdfstartview=FitH,
pdfhighlight=/N
}

\makeatletter
\def\@captionheadfont{\scshape}
\def\@captionfont{\normalfont\small}
\makeatother

\newcommand{\inclgr}{\includegraphics[scale=1.5]}
\newcommand{\inclgraph}{\includegraphics[scale=1.1]}

\theoremstyle{plain}
\newtheorem{conj}{Conjecture}
\newtheorem{theoremSta}{Stahl Theorem}
\newtheorem{theoremBusTwoPoint}{Buslaev Two-Point Theorem}

\theoremstyle{definition}
\newtheorem{remark}{Remark}
\newtheorem{Case}{Case}

\def\({\left(}
\def\){\right)}
\def\GRS{\operatorname{GRS}}
\def\mcap{\operatorname{cap}}
\def\supp{\operatorname{supp}}
\def\const{\operatorname{const}}
\def\bus{\operatorname{Bus}}
\def\NN{\mathbb N}
\def\RR{\mathbb R}
\def\CC{\mathbb C}
\def\QQ{\mathbb Q}
\def\PP{\mathbb P}
\def\ZZ{\mathbb Z}
\def\UU{U}
\def\FF{\mathscr T}
\def\GG{\mathscr G}
\def\HH{\mathscr H}
\def\myA{\mathscr A}
\def\RS{\mathfrak R}
\def\myf{\mathfrak f}
\def\zz{\mathbf z}
\let\leq\leqslant
\let\geq\geqslant

\let\myh\widehat
\let\myo\overline

\begin{document}

\title{On the limit zero distribution of type~I Hermite--Pad\'e polynomials}

\author{Nikolay~R.~Ikonomov}
\address{Institute of Mathematics and Informatics, Bulgarian Academy of
Sciences}
\email{nikonomov@math.bas.bg}
\author{Ralitza~K.~Kovacheva}
\address{Institute of Mathematics and Informatics, Bulgarian Academy of
Sciences}
\email{rkovach@math.bas.bg}
\author{Sergey~P.~Suetin}
\address{Steklov Mathematical Institute of Russian Academy of Sciences, Russia}
\email{suetin@mi.ras.ru}
\thanks{The third author was partially supported by the Russian Foundation for Basic Research
(RFBR, grants \textnumero 13-01-12430-ofi-m2 and \textnumero 15-01-07531-a),
and Russian Federation Presidential Program for support of the leading scientific schools (grant
NSh-2900.2014.1).}

\date{June 26, 2015}

\begin{abstract}
In this paper are discussed the results of new numerical experiments
on zero distribution of type~I Hermite--Pad\'e polynomials of order $n=200$
for three different collections of three functions $[1,f_1,f_2]$.
These results are obtained by the authors numerically
and do not match any of the theoretical results that were proven
so far. We consider three simple cases of multivalued analytic
functions $f_1$ and $f_2$, with separated pairs of branch points
belonging to the real line. In the first case both functions have
two logarithmic branch points, in the second case they both have
branch points of second order, and finally, in the third case they
both have branch points of third order.

All three cases may be considered as representative of the
asymptotic theory of Hermite--Pad\'e polynomials. In the first two
cases the numerical zero distribution of type~I
Hermite--Pad\'e polynomials are similar to each other, despite the
different kind of branching. But neither the logarithmic case,
nor the square root case can be explained from the asymptotic
point of view of the theory of type~I Hermite--Pad\'e polynomials.

The numerical results of the current paper might be considered as
a challenge for the community of all experts on Hermite--Pad\'e
polynomials theory.

Bibliography: \cite{Wid69}~items.

Figures: \ref{Fig_cub(p8s)2000(200)_bk}~items.
\end{abstract}

\maketitle


\markright{ZERO DISTRIBUTION OF HERMITE--PAD\'E POLYNOMIALS}

\setcounter{tocdepth}{2}
\tableofcontents

%%%fulltext

\section{Introduction}\label{s1}

In this article, we are concerned with the general problem of the
limit zero distribution (LZD) of type~I Hermite--Pad\'e (HP)
polynomials for the collection of $m$ function
$[f_0\equiv1,f_1,\dots,f_{m-1}]$, $m\geq3$, where
$f_1,\dots,f_{m-1}$ are multivalued analytic functions, each one
having a finite set of branch points on the Riemann sphere
$\myo\CC$; the reader is referred at first to the basic
surveys~\cite{Nut84}, \cite{Sta88}, and
also~\cite{Gon03}, \cite{NiSo88}, \cite{Apt08}. Up to the end of the paper,
suppose that the functions $f_0\equiv1,f_1,\dots,f_{m-1}$ are
rationally independent (over the field $\CC(z)$).

We restrict our attention to the case of the collection of three
functions, in other words, to the case $m=3$. Let us recall some
basic notions from the HP polynomials theory
(see~\cite{Nut84},~\cite{Sta88}).

Suppose that two functions $f_1$ and $f_2$ are given by convergent
power series at the infinity point $z=\infty$. For a fixed number
$n\in\NN$, the type~I HP polynomials $Q_{n,0},Q_{n,1},Q_{n,2}$ of
degree $\leq{n}$, $Q_{n,j}=Q_{n,j}(z;f_1,f_2)\in\PP_n$,
$\PP_n:=\CC_n(z)$, $j=0,1,2$, not all $Q_{n,j}\equiv0$, are
defined in such a way that the following basic relation holds
true:
\begin{equation}
R_n(z):=(Q_{n,0}\cdot1+Q_{n,1}f_1+Q_{n,2}f_2)(z)=O\(\frac1{z^{2n+2}}\),
\quad z\to\infty.
\label{1.1}
\end{equation}
It is well known that such polynomials $Q_{n,j}$ always exist and
in the ``generic case''\footnote{In other terminology, ``in common position''.} they
are unique up to a suitable normalization
(see~\cite{Nut84}, \cite{Sta88}, \cite{NiSo88}, \cite{FiLo11}).

To be more precise, we treat here, from numerical point of view, three different
pairs of Markov-type functions: $f_1,f_2$, $g_1,g_2$ and $h_1,h_2$.

Recall that a function $f\in\HH(\infty)$ is a Markov-type function if
$f(z)=\myh{\mu}(z)+\const$, where $\myh{\mu}$ is the
Cauchy transform of a positive Borel measure, with support $\supp\mu\Subset\RR$
(a compact subset of the real line $\RR$),
\begin{equation}
\myh{\mu}(z):=\int\frac{d\mu(x)}{z-x},\quad z\not\in\supp\mu.
\label{mar1}
\end{equation}
Up to the end of the paper, a support of a Markov-type function
$f=\myh{\mu}+\const$ means the support of the corresponding
measure~$\mu$.

Let $E_1:=[a_1,a_2]=[-1,a]$, $E_2:=[b_1,b_2]=[-a,1]$, where
$a\in(0,1)$ is a real parameter. Hence both segments $E_1$ and
$E_2$ are overlapping, $E_1\cap E_2=[-a,a]\not=\varnothing$.

\begin{Case}\label{cas1}
Let
\begin{align}
f_1(z)&:=\int_{-1}^a\frac{dx}{z-x}=\log\frac{z-a}{z+1},\quad z\notin E_1,
\label{f1}\\
f_2(z)&:=\int_{-a}^1\frac{dx}{z-x}=\log\frac{z-1}{z+a},\quad z\notin E_2.
\label{f2}
\end{align}
We take the main branch of the logarithmic function in~\eqref{f1} and~\eqref{f2},
in the sense that
$$
\log\frac{z-a}{z+1},\ \log\frac{z-1}{z+a}\approx\log{1}=0,\quad\text{as}\quad
z\to\infty.
$$
\end{Case}

\begin{Case}\label{cas2}
Let
\begin{align}
g_1(z)
&:=\(\frac{z-a}{z+1}\)^{1/2}=\frac1\pi\int_{-1}^a\sqrt{\frac{a-x}{x+1}}
\frac{dx}{x-z}+1,\quad z\notin E_1,
\label{1.3}\\
g_2(z)
&:=\(\frac{z-1}{z+a}\)^{1/2}=\frac1\pi\int_{-a}^1\sqrt{\frac{1-x}{x+a}}
\frac{dx}{x-z}+1,\quad z\notin E_2.
\label{1.4}
\end{align}
We take such a branch of the $(\cdot)^{1/2}$
function, that $g_1(z),g_2(z)\to1$ as $z\to\infty$, and under
the square root function $\sqrt{\,\cdot\,}$ we mean the ``arithmetic
square root function'', that is $\sqrt{x^2}=x$ for $x\in\RR_{+}$.
\end{Case}

\begin{Case}\label{cas3}
Let
\begin{align}
h_1(z)&:=\(\frac{z-a}{z+1}\)^{1/3},\quad z\notin E_1,
\label{c3.1}\\
h_2(z)&:=\(\frac{z-1}{z+a}\)^{1/3},\quad z\notin E_2.
\label{c3.2}
\end{align}
We take such a branch of the $(\cdot)^{1/3}$
function, that $h_1(z),h_2(z)\to1$ as $z\to\infty$, and, in what follows,
under the cubic root function $\sqrt[3]{\,\cdot\,}$ we mean the ``arithmetic
cubic root function'', that is $\sqrt[3]{x^3}=x$ for $x\in\RR_{+}$.
\end{Case}

Clearly, all three functions $f_j,g_j,h_j$, $j=1,2,3$ are
Markov-type functions. But since $a\in(0,1)$, we have $E_1\cap
E_2=[-a,a]\neq\varnothing$ and $E_1\neq E_2$. Hence, all three
pairs $f_1,f_2$, $g_1,g_2$ and $h_1,h_2$ are neither an Angelesco,
nor a Nikishin system (see~\cite{Nut84},~\cite{GoRa81},
\cite{NiSo88},~\cite{GoRaSo97}, \cite{Gon03},~\cite{ApLy10}). Thus
so far, not even one theorem of general type is known, which
predicts the LZD of type~I HP polynomials for the three
collections of multivalued analytic functions $[1,f_1,f_2]$,
$[1,g_1,g_2]$ and $[1,h_1,h_2]$, that were introduced above.
Hence, the numerical results that are presented in this paper
should be explained from theoretical point of view by the
forthcoming papers.

Since $a$ is a parameter, we actually have three families of
Markov-type functions: $\FF=\{f_1=f_1(z;a),\ f_2=f_2(z;a),\
a\in(0,1)\}$, $\GG=\{g_1=g_1(z;a),\ g_2=g_2(z;a),\ a\in(0,1)\}$,
and $\HH=\{h_1=h_1(z;a),\ h_2=h_2(z;a),\ a\in(0,1)\}$. Thus, for a
fixed $n\in\NN$, we have three families of HP polynomials defined
by the basic relation~\eqref{1.1}. To distinguish between
Case~\ref{cas1}, Case~\ref{cas2}, and Case~\ref{cas3}, we
introduce new notations for the HP polynomials for the collections
of three functions $[1,g_1,g_2]$ and $[1,h_1,h_2]$. We set %, namely,
$P_{n,0},P_{n,1},P_{n,2}\in\PP_n\setminus\{0\}$ for the
HP polynomials corresponding to the collection $[1,g_1,g_2]$,
\begin{equation}
(P_{n,0}\cdot1+P_{n,1}g_1+P_{n,2}g_2)(z)=O\(\frac1{z^{2n+2}}\),
\quad z\to\infty,
\label{1.5}
\end{equation}
and $\UU_{n,0},\UU_{n,1},\UU_{n,2}\in\PP_n\setminus\{0\}$ for the
HP polynomials corresponding to the collection $[1,h_1,h_2]$,
\begin{equation}
(\UU_{n,0}\cdot1+\UU_{n,1}h_1+\UU_{n,2}h_2)(z)=O\(\frac1{z^{2n+2}}\),
\quad z\to\infty.
\label{1.5.2}
\end{equation}

Notice that for each value of the parameter $a\in(0,1)$ the
supports of the two Markov-type functions $f_j$ and $g_j$,
$j=1,2$, coincide with each other, but the types of the
corresponding branch points are of different nature. It is worth
noting that in each pair $f_1,f_2$ and $g_1,g_2$, the branch
points are separated from each other
(cf.~\cite{NuTr87},~\cite{ApKuAs07}). But since for each
$a\in(0,1)$, $E_1\cap E_2\neq\varnothing$, $E_1\not\subset E_2$,
and $E_2\not\subset E_1$, the system $f_1,f_2$ is neither an
Angelesco, nor a Nikishin system. The same is true for the
systems $g_1,g_2$ and $h_1,h_2$.

In all of the numerical experiments discussed in the present paper, we
set $n=200$ (see~\eqref{1.1}).
%For the numerical calculations we
%set $n=200$ (see~\eqref{1.1}) and follow that up to the end of the paper.
We also fix the set of
values of parameter $a$ as to be
$a=0.2;0.4;0.625;0.73;0.8$.\footnote{For
Case~3 we set $a=0.85$ instead of $0.8$.}

At first we set $a=-0.1$ (and keep $n=200$). Clearly, that value
of~$a$ is out of the range of $(0,1)$. But the reason is that for
$a=-0.1$ all three pairs $f_1,f_2$, $g_1,g_2$ and $h_1,h_2$ form
Angelesco systems, since $E_1,E_2\subset\RR$ and $E_1\cap
E_2=\varnothing$. On Fig.~\ref{Fig_log(pL1s)1500(200)_full} and
Fig.~\ref{Fig_sqrt(pL1s)2000(200)_full} zeros of HP polynomials
$Q_{200,j}$ and $P_{200,j}$, $j=0,1,2$, are plotted. The numerical
distribution of the zeros of the HP polynomials are similar in
both cases and do not depend on the different type of branching of
the functions $f_1,f_2$ and $g_1,g_2$ at the points $z=\pm1,\pm
a$. The reason is as follows. All these numerical results are in
full agreement with the general theory of limit zero distribution
(LZD) for HP polynomials in the Angelesco case (see~\cite{GoRa81},
\cite{GoRaSo97},~\cite{ApLy10}; in these papers the case of
type~II HP polynomials was considered, but it is well known that
in the Angelesco's case there is a direct connection between LZD
of type~I and type~II HP polynomials). Since the sizes of the
supports of Markov-type functions $f_1$ and $f_2$ (and $g_1$,
$g_2$ as well) are equal, the phenomena of the so-called ``pushing
of the charge''\footnote{The phenomena was discovered in 1981 by
Gonchar and Rakhmanov in~\cite{GoRa81}.} is absent in the case
$a=-0.1$. It is in full agreement also with the general theory
(see~\cite{Kal79},~\cite{ApKa86},~\cite{GoRa81}). From
Fig.~\ref{Fig_log(pL1s)1500(200)_full},
Fig.~\ref{Fig_log(pL1s)1500(200)_blue},
Fig.~\ref{Fig_sqrt(pL1s)2000(200)_full},
Fig.~\ref{Fig_sqrt(pL1s)2000(200)_blue} follows, that the zeros of
HP polynomials $Q_{200,1}$ and $Q_{200,2}$ (and $P_{200,1}$,
$P_{200,2}$ as well) are located on the segments $E_1$ and $E_2$,
respectively. But the zeros of $Q_{200,0}$ (and $P_{200,0}$
respectively) are located on the imaginary axis. Thus in some
sense the zeros of $Q_{200,0}$ ``intend to separate'' the zeros of
$Q_{200,1}$ and $Q_{200,2}$ from each other. Again this result is
in good agreement with the well known LZD of type~I HP polynomials
in the Angelesco case (see, e.g.~\cite{NuTr87}). To finish the
analysis of the case $a=-0.1$, we emphasize that in this case all
numerical results are in good agreement with the general theory of
HP polynomials. Thus, it is reasonable to consider all numerical
results for $a=0.2;0.4;0.625;0.73;0.8$ as to be trustable.

We describe in short the numerical experiments.

We set the parameter $a$ to be $0.2;0.4;0.625;0.73;0.8$, and for
each value of $a$ we find numerically the zeros of HP polynomials
$Q_{200,j}(z;a)$, $P_{200,j}(z;a)$ and $\UU_{200,j}(z;a)$,
$j=0,1,2$, respectively. Thus we obtain a ``series'' of numerical
experiments for the different values of the parameter $a$. By the
reasons mentioned above, it is meaningful to consider all these
numerical results as trustable, and hence we can analyze them
from a theoretical point of view.

Before laying out in details the discussion and the analysis of
the results of the numerical experiments (see
subsection~\S\ref{s4s3}), we recall basic general facts about
Pad\'e approximants\footnote{Sometimes in the current paper we
follow the Chudnovsky's terminology~\cite{ChCh84}: under ``Pad\'e
approximants'' (and ``Pad\'e polynomials'' respectively), we also
mean the Hermite--Pad\'e approximants, the multipoint Pad\'e
approximants, etc.} convergence theory (PA-theory). For this purpose, we provide in the
next sections~\S\ref{s2} and~\S\ref{s3} a survey for
the most important results of PA convergence theory, which is intended to
make more clear the mathematical essence of the numerical results
presented in the current paper.

\section{Classical Pad\'e approximants (or \texorpdfstring{$J$}{J}-fractions) convergence theory}\label{s2}

\subsection{Stahl's \texorpdfstring{$S$}{S}-curve and Stahl's limit zero-pole
distribution theorem in classical Pad\'e approximants theory}\label{s2s1}
%\subsection{Classical Pad\'e approximants: Stahl's \texorpdfstring{$S$}{S}-curve and Stahl's limit zero-pole
%distribution theorem}\label{s2s1}

The basic general results of PA-theory are of two types. The first type is
concerned with the limit zero distribution of Pad\'e polynomials for
multivalued analytic functions with a finite set of branch points on
the Riemann sphere $\myo\CC$. Commonly, these results are referred
to as ``weak asymptotics of PA''.
The second type of results are devoted to the so-called ``strong asymptotics of PA''.

The problem of first type
was originally stated in a formal way by J.~Nuttall (see~\cite{Nut80}
and the review~\cite{Nut84}) at the beginning of 1970's, and for the
classical (``one-point'') Pad\'e approximants the problem was completely
solved\footnote{It is worth noting that seminal Stahl Theorem is
much more general than the original conjecture of Nuttall.
Stahl Theorem was proved under the condition that the set of singular
points of the given multivalued function is of zero (logarithmic) capacity
instead of being a finite set.} by
H.~Stahl~\cite{Sta85a}--\cite{Sta86b} in the middle of 1980's (see
also~\cite{Sta97b}).

The second type problem was also initiated by Nuttall
(see~\cite{NuSi77}, \cite{Nut84}, \cite{Nut86}, \cite{Nut90} and
also~\cite[\S~2, Conjecture~1]{KovSu14}). Recently were obtained
strong results in this direction
(see~\cite{Sue00},~\cite{ApYa11},~\cite{MaRaSu12}). However,
generally the problem is still far from being solved even for the
classical PA. We explain the main reason for this situation.
The ``strong asymptotics theory'' is a
far and deep generalization of the classical
Bernstein--Szeg\H{o} asymptotic theory for the ordinary
polynomials, orthogonal on the unit segment $\Delta:=[-1,1]$ and
related to the general case of non-Hermitian orthogonal
polynomials (see~\cite{Sze62},~\cite{Akh60},~\cite[Section~1 and
Section~6]{Nut84},~\cite{Nut90}, the paper by
H.~Widom~\cite{Wid69}, as well).

Given a germ\footnote{For example, a convergent power series at a point
$z_0\in\myo\CC$; commonly we set $z_0=\infty$ or $z_0=0$.} $f$
of a multivalued analytic function $f$ with a finite number of
branch points, the seminal Stahl Theorem\footnote{Because of its
general character and the various subjects of complex analysis
involved into the proof of the theorem, it is sometimes considered
as ``Stahl Theory''.} gives a complete answer to the problem of
limit zero-pole distribution of the classical PA to~$f$. Stahl
Theorem is of a quite general character. It not only admits
multivalued functions with a finite number of branch points on the
Riemann sphere, but also multivalued analytic functions with a
singular set of zero (logarithmic) capacity. The keystone of Stahl
theory~\cite{Sta85a}--\cite{Sta86b} is the existence of a unique
(up to a compact set of zero capacity) maximal domain for the given
at the point $z = \infty$ multivalued function $f$. The so-called ``maximal domain'' of
holomorphy of $f$ is a domain $D=D(f)\ni\infty$, such that the
given germ $f$ can be continued as holomorphic (analytic and
single-valued) function from a neighborhood of the infinity point
$z=\infty$ to $D$ (the function~$f$ can be continued analytically
along each path that belongs in to $D$). ``Maximal'' means that
$\partial{D}$ is of ``minimal capacity'' among all compact sets
$\partial{G}$, such that $G$ is a domain, $\infty \in G$ and
$f\in\HH(G)$. Such ``maximal'' domain $D$ is unique (up to a
compact set of zero capacity). The compact set
$S=S(f):=\partial{D}$ is now called ``Stahl's compact set'' or
``Stahl's $S$-compact set'', and $D$ is called ``Stahl's domain''.
The crucial properties of $S$ for the Stahl theory to be true
are the following: the complement $D=\myo\CC\setminus{S}$ is a
domain, $S$ consists of a finite number of analytic arcs (in fact,
the union of the closures of the critical trajectories of a
quadratic differential), and, finally, $S$ possesses the so-called
property of ``symmetry''\footnote{Compact sets of such type are
usually referred to as ``$S$-compact sets'' or ``$S$-curves'',
see~\cite{MaRa11}, \cite{RaSu13},~\cite{Rak12},~\cite{KoIkSu15}.}, that is,
\begin{equation}
\frac{\partial g^{}_D(z,\infty)}{\partial n^{+}}
=\frac{\partial g^{}_D(z,\infty)}{\partial n^{-}},\quad z\in S^0,
\label{1.6}
\end{equation}
where $g_D(z,\infty)$ is Green's function for the domain $D$,
with the logarithmic singularity at the point $z=\infty$, $S^0$ is the
union of all open arcs of~$S$ (which closures constitute $S$, that is
$S\setminus S^0$ is a finite\footnote{In the general Stahl Theorem the
set $S\setminus S^0$ is of zero capacity.} set), and $\partial n^{+}$
and $\partial n^{-}$ are the inner (with respect to $D$) normal
derivatives of $g_D(z,\infty)$ at a point $z\in S^0$ from the opposite
sides of $S^0$.

Given a finite set $\Sigma\subset\myo\CC$,
$\#\Sigma<\infty$, let $\myA(\myo\CC\setminus\Sigma)$
be the set of all functions analytic in the domain
$\myo\CC\setminus\Sigma$. Let
$\myA^0(\myo\CC\setminus\Sigma):=\myA(\myo\CC\setminus\Sigma)
\setminus\HH(\myo\CC\setminus\Sigma)$, i.e. for a fixed $\Sigma$ a
function $f$ is from the set $\myA^0(\myo\CC\setminus\Sigma)$ if it is a
multivalued analytic function in the domain
$\myo\CC\setminus\Sigma$, but not a holomorphic function in
$\myo\CC\setminus\Sigma$ ($f$ is analytic, but not single-valued).

Let $\Sigma=\Sigma_f=\{a_1,\dots,a_p\}$, $\#\Sigma=p<\infty$, be the set of
all branch points of~$f$, i.e. $f\in\myA(\myo\CC\setminus\Sigma)
\setminus\HH(\myo\CC\setminus\Sigma)=:\myA^0(\myo\CC\setminus\Sigma)$. Clearly,
if $\Sigma_f=\{-1,1\}$, then Stahl's compact set (with respect to the infinity
point) $S = [-1,1]$.

The properties of the compact set $S$ described above are crucial
for the validness of Stahl Theorem.

We recall some basic notations from the general PA-theory.

Let there be given a function $f\in\myA^0(\myo\CC\setminus\Sigma)$, which is holomorphic
at the infinity point, $f\in\HH(\infty)$. To be more precise, suppose we
are given a convergent power series, i.e. a germ,
\begin{equation}
f=\sum_{k=0}^\infty\frac{c_k}{z^{k+1}} ,
\label{1.7}
\end{equation}
and the corresponding analytic function
$f\in\myA^0(\myo\CC\setminus\Sigma)$, where $\#\Sigma<\infty$. For a
positive Borel measure $\mu$ with a compact support
$\supp(\mu)\Subset\myo\CC$, denote by $V^\mu(z)$ the logarithmic
potential\footnote{If necessary, the potential should be spherically
normalized, see~\cite{GoRa87},~\cite{Bus13} and Remark~\ref{rem7}.} of
$\mu$, that is:
$$
V^\mu(z):=\int_{\supp\mu}\log\frac1{|z-\zeta|}\,d\mu(\zeta).
$$

Given a number $n\in\NN$, let\footnote{There are some ambiguities
in this notations with~\eqref{1.5}, apologies for that.}
$P_{n,0},P_{n,1}\in\CC_n(z)$, $P_{n,1}\not\equiv0$, be the Pad\'e polynomials (at
the infinity point) of the function $f$. We recall that
\begin{equation}
R_n(z):=(P_{n,0}+P_{n,1}f)(z)=O\(\frac1{z^{n+1}}\),\quad z\to\infty.
\label{1.8}
\end{equation}
The polynomials $P_{n,0},P_{n,1}$ are not unique, but the ratio
$P_{n,0}/P_{n,1}$ is uniquely determined by the relation~\eqref{1.8}.
The rational function $[n/n]_f:=-P_{n,0}/P_{n,1}$ is called the
diagonal Pad\'e approximant of order $n$ (or the $n$-diagonal Pad\'e
approximant) of the function~$f$ (at the infinity point),
and $R_n$ is the remainder function of the Pad\'e approximant.

Given an arbitrary polynomial $Q\in\CC(z)$, $Q\not\equiv0$, denote by
$$
\chi(Q):=\sum_{\zeta:Q(\zeta)=0}\delta_\zeta
$$
the associated zero counting measure of $Q$ (as
usual, $\delta_\zeta$ denotes the Dirac measure concentrated at the point
$\zeta\in\myo\CC$).

One of the main results of Stahl theory~\cite{Sta85a}--\cite{Sta86b} is

\begin{theoremSta}[H.~Stahl, see~\cite{Sta86a},~\cite{Sta97b}]
Let $f\in\HH(\infty)$, $f\in\myA^0(\myo\CC\setminus\Sigma$),
$\#\Sigma<\infty$. Let $D=D(f)$ be Stahl's ``maximal'' domain for $f$,
$S=\partial{D}$ -- Stahl's compact set,
$[n/n]_{f}=-P_{n,0}/P_{n,1}$ -- the $n$-diagonal Pad\'e approximant to
the function~$f$. Then the following statements are valid:

1) there exists LZD of Pad\'e polynomials $P_{n,j}$,
$j=0,1$, namely,
\begin{equation}
\frac1n\chi(P_{n,j})
\overset{*}\longrightarrow\lambda,\quad\text{as}\quad
n\to\infty,\quad j=0,1,
\label{2.1}
\end{equation}
where $\lambda=\lambda_S$ is a unique probability equilibrium
measure of the compact set~$S$;
%that is $V^\lambda(z)\equiv\gamma_S$, $z\in S$,
%where $\gamma_S$ is the Robin constant for~$S$.

2) the $n$-diagonal Pad\'e approximants converge in capacity to the
function~$f$ inside
%\footnote{That is, on compact subsets of the domain~$D$.}
the domain~$D$,
\begin{equation}
[n/n]_f(z)\overset{\mcap}\longrightarrow f(z),\quad n\to\infty,\quad z\in D;
\label{2.2}
\end{equation}

3) the rate of the convergence in~\eqref{2.2} is completely characterized by the equality
\begin{equation}
\bigl|(f-[n/n]_f)(z)\bigr|^{1/n}\overset{\mcap}\longrightarrow
e^{-2g^{}_D(z,\infty)},\quad n\to\infty,\quad z\in D.
\label{2.3}
\end{equation}
\end{theoremSta}

\begin{remark}\label{rem1}
In item~1), $\lambda_S$ is the unique probability equilibrium
measure of the compact set $S$, that is $V^{\lambda_S}(z)\equiv\gamma_S,\, z\in S$,
where $V^{\lambda_S}(z)=-\int\log|z-\zeta|\,d\lambda_S(\zeta)$ is the
logarithmic potential of the measure $\lambda_S$, $\gamma_S$ is the
Robin constant for~$S$. The notation ``$\overset{*}\longrightarrow$''
stands for convergence of measures in the weak-star topology. The notation
``$\overset{\mcap}\longrightarrow$'' in items~2) and~3) means
convergence in capacity inside (on compact subsets of) the domain~$D$.

As it is basely known,
$g_D(z,\infty)\equiv\gamma_S-V^{\lambda_S}(z)$. From here it follows
immediately, that the $S$-property~\eqref{1.6} may be written
in an equivalent form (cf.~\eqref{b3}), namely
\begin{equation}
\label{1.6eq}
\frac{\partial V^\lambda}{\partial n^{+}}(z)
=\frac{\partial V^\lambda}{\partial n^{-}}(z),\quad z\in S^0.
\end{equation}
\end{remark}

\begin{remark}\label{rem2}
Let $P^*_{n,j}$, $j=0,1$, be the monic Pad\'e polynomials.
Then the statement of item~1) is equivalent to the following:
\begin{equation}
\bigl|P^*_{n,j}(z)\bigr|^{1/n}\overset{\mcap}\longrightarrow
e^{-V^{\lambda_S}(z)
},\quad n\to\infty,\quad z\in D, \quad j=0,1.
\label{2.4}
\end{equation}
Similarly, denote by $R^{*}_n(z)$ the normalized remainder
function,
$$
R^{*}_n(z)=\frac1{z^{n+1+\ell_n}}+\dotsb,\quad z\to\infty,
$$
where $\ell_n\in\ZZ$, $\ell_n=o(n)$ as $n\to\infty$.
Then the statement of item~3) is equivalent to the following:
\begin{equation}
\bigl|R^{*}_n(z)\bigr|^{1/n}\overset{\mcap}\longrightarrow
e^{-g^{}_D(z,\infty)},\quad n\to\infty,\quad z\in D.
\label{2.5}
\end{equation}
\end{remark}

Let $\Sigma_f=\{-1,1\}$. Clearly, in this $S = [-1,1]$.
Hence, for the case of a general multivalued analytic function
$f\in\myA^0(\myo\CC\setminus\Sigma_f)$, $\#\Sigma_f<\infty$, the
associated (with the function $f$) Stahl's compact set $S$ should
be considered as a natural replacement of the unit segment
$\Delta=[-1,1]$, which is the main object in the classical
Bernstein--Szeg\H{o} asymptotic theory. Recall that the classical
Bernstein--Szeg\H{o} formula of strong asymptotics for
polynomials, orthogonal on the segment $\Delta$ with respect to a weight $\rho$ defined on $\Delta$
and ``smooth'' enough, is
presented in terms of Szeg\H{o} function, associated with
$\rho$, and the ``mapping'' function $\varphi(z)=z+(z^2-1)^{1/2}$,
associated with $\Delta$. Viewing
Stahl Theorem, it is natural to conjecture, that in the formula of strong asymptotics of
Pad\'e polynomials $P_{n,j}$ Stahl's compact set $S$ should appear
in a natural way, but not the segment $\Delta$
or the finite union of the real segments
(see~\cite{Akh60},~\cite{Wid69},~\cite{Nut84},~\cite[Conjecture~1]{Nut86},
\cite{ApYa11},~\cite[Conjecture~2]{KovSu14}).

The existence of Stahl's $S$-curve for an arbitrary multivalued
analytic function with a finite set of singular points, on one hand, and Stahl
Theorem of the LZD of Pad\'e polynomials, on the other, should play a crucial
role in the forthcoming development of the theory of the {\it
strong asymptotics} of Pad\'e polynomials, associated with a
multivalued analytic function, and the corresponding theory of
{\it strong convergence} of Pad\'e approximants
(cf.~\cite{GoSu04},~\cite{ApYa11}).

It is worth noting that, in general, the property of a symmetry
(in other words, $S$-property of type~\eqref{1.6},~\eqref{b3},
see~\cite[formula~(11)]{RaSu13}) gives rise to a local analogue of
complex conjugation operation and to the method of the
construction of an associated two-sheeted or three-sheeted
Riemann surface
(see~\cite{Sta97b},~\cite{GoRaSu91},~\cite[\S~6]{RaSu13}).

In 1987 A.~A.~Gonchar and E.~A.~Rakhmanov~\cite{GoRa87} extended
the original Stahl's conception of symmetry (the $S$-curve
conception) to the conception of ``symmetry in the presence of an
external field''.\footnote{See also the pioneering paper by
Gonchar and Rakhmanov~\cite{GoRa81}, 1981, where for the first time the
notion of logarithmic theoretical potential problem with an
external field was introduced in connection
with the problem of LZD of type~II HP polynomials for the
Angelesco systems of Markov-type functions.} The general LZD-theorem of
orthogonal polynomials was proven mainly by using the new notions in~\cite{GoRa87};
i.e., for polynomials orthogonal with respect
to a variable weight (depending on the degree of the polynomial, see
also~\cite{GoRa84}). This new concept of the ``presence of an external
field'' offered Gonchar and Rakhmanov the chance to solve the
well known ``$1/9$-problem'' (see~\cite{Var74}, also
\href{http://mathworld.wolfram.com/One-NinthConstant.html}{Wolfram MathWorld: One-Ninth Constant},
it is also known as ``Varga's constant'', see
\href{http://mathworld.wolfram.com/VargasConstant.html}{Wolfram MathWorld: Varga's Constant}).
Finally, Stahl's and Gonchar--Rakhmanov's results gave rise to the
so-called $\GRS$-method (Gonchar--Rakhmanov--Stahl-method), which solves
the problem of LZD of generalized (non-Hermitian) quasi-orthogonal polynomials.
The method itself consists of three steps.

First step: to state an appropriate theoretical potential
extremal problem associated with the initial problem of LZD of
Pad\'e polynomials\footnote{That is, Pad\'e polynomials,
multipoint Pad\'e polynomials, Hermite--Pad\'e polynomials,
generalized Pad\'e polynomials, etc., see~\cite{ChCh84}.}, and
then proving the existence of an ``extremal'' compact set in some
family of ``admissible'' compact sets.\footnote{Sometimes it is
called a ``stationary'' compact set of the associated extremal
theoretic potential problem, see~\cite{MaRaSu11b},~\cite{RaSu13}.}

Second step: to prove that the extremal compact set possesses
the special property of ``symmetry'', i.e., the set is a
``weighted'' $S$-curve, associated with the problem under consideration.

Third step: to establish that the LZD of the generalized orthogonal
polynomials under question exists and coincides with the equilibrium
 probability measure (in the presence of the external field), which is
concentrated on the weighted $S$-curve.

\subsection{Strong asymptotics in classical Pad\'e approximants theory}\label{s2s2}

% Here
We give one of the possible forms of the strong asymptotics
of Pad\'e polynomials under the conditions of Stahl Theorem
(cf.~\cite[\S~2, Conjecture~1]{KovSu14}). It is worth noting that
in general (except for the case of genus zero) such a
representation is not unique
(see~\cite{DeKrMc99b},~\cite{Pas06},~\cite{KoSu14b}). The reason
is that there exist so-called spurious zeros of Pad\'e
polynomials, with behavior, as $n\to\infty$, that can be described
in many ways.

\begin{conj}[cf.~\cite{KovSu14}, \rm\S~2, Conjecture~1]\label{conj1}
Let $f\in\HH(\infty)$ and $f\in\myA^0(\myo\CC\setminus\Sigma)$ for some
finite set $\Sigma\subset\myo\CC$. Let $D=D(f)$ be Stahl's maximal
domain associated with $f$, $S=S(f)=\partial{D}$ -- the
corresponding Stahl's compact set for~$f$, $\RS_2=\RS_2(f)$ --
the canonical\footnote{In Stahl's sense, see~\cite{Sta97b} and
also~\cite{ApBuMaSu11}, \cite{KovSu14}.} hyperelliptic two-sheeted
Riemann surface associated\footnote{Such Riemann surface is uniquely
determined by~$S$.} with the compact set~$S$,
$\zz=z^{(1,2)}=(z,\pm)\in\RS_2$ -- an arbitrary point on the
two-sheeted $\RS_2$, $\Psi_n(\zz)=\Psi_n(\zz;f)$ -- Nuttall's
psi-function associated\footnote{See~\cite{Sue00},~\cite{Sue06},
\cite{ApBuMaSu11},~\cite{KovSu14}.} with $f$ and $\RS_2$. Then, after a
suitable normalization of the Pad\'e polynomials $P_{n,j}(z)=P_{n,j}(z;f)$,
$j=0,1$, and the remainder function $R_n$, the following relations take
place in capacity inside the domain~$D$:
\begin{align}
P_{n,j}(z)&\overset{\mcap}=\frac{(-1)^j}{f^j(z)}\Psi_n(z^{(1)})\bigl(1+o(1)\bigr),\quad &n\to\infty,
\nonumber \\
R_n(z)&\overset{\mcap}=\frac{\Psi_n(z^{(2)})}{w(z^{(2)})}
\bigl(1+o(1)\bigr),\quad &n\to\infty,
\label{2.6}
\end{align}
where $w^2=H_{2g+2}(z)$, $H_{2g+2}\in\CC_{2g+2}[z]$, is the equation that
determines\footnote{Recall that all zeros of the polynomial
$H_{2g+2}$ are simple.} the hyperelliptic Riemann surface $\RS_2$ of genus~$g$.
\end{conj}

Conjecture~\ref{conj1} is based partly on Nuttall's result from
1986 \cite{Nut86} (see also \cite{MaRaSu12})
via Liouville--Steklov asymptotic method for
the function $f$ of the special form
\begin{equation}
f(z)=\prod_{j=1}^3(z-a_j)^{\alpha_j},\quad\sum_{j=1}^3\alpha_j=0,
\quad\alpha_j\in\CC\setminus\ZZ,\quad f(\infty)=1
\label{2.7}
\end{equation}
(cf.~\cite[\S~5. Asymptotic Conjecture]{Nut86}).

Conjecture~\ref{conj1} was proven in some particular cases
in~\cite{Sta96}, \cite{Sue00}, \cite{ApYa11},
comp. also~\cite{Nut84}, \cite{Nut90}, \cite{ApBuMaSu11}, \cite{KovSu14}.

In general, the problem of the strong asymptotics of
Pad\'e polynomials is still open and the state is as follows:
there are only partial results about the subject,
and all of them are based on the existence of
Stahl's compact set $S$ and the concepts associated with it, namely, the
two-sheeted hyperelliptic Riemann surface and the corresponding Abelian
integrals on it. The problem does not depend
on the method used to produce the asymptotic formula:
method of matrix Riemann--Hilbert boundary value
problem, Liouville--Steklov method, method of the Nuttall's singular
integral equation. Ultimately, all these methods are based on
the existence of Stahl's $S$-curve.

\section{Two-point Pad\'e approximants (or \texorpdfstring{$T$}{T}-fractions) convergence theory}\label{s3}

Recently V.~I.~Buslaev~\cite{Bus13} (see also~\cite{BuMaSu12}, \cite{Bus15})
applied the $\GRS$-method in treating the multipoint\footnote{To be
more precise, $m$-point where $m\in\ZZ_{>0}$ is fixed.} PA of a
multivalued analytic function with a finite set of branch points.
Buslaev~\cite{Bus13} proved an analogue of the classical Stahl
Theorem, furthermore, he discovered some special features of
the multipoint PA. In particular, for the case of two-point PA, say at
the points $z=0$ and $z=\infty$, he was the first, who found
that in the ``generic case'' the weighted $S$-curve\footnote{It is
natural to call such a curve a ``Buslaev's $S$-compact set''.},
associated with this problem, partitions the Riemann sphere in an
``optimal way'' into two domains $D_0\ni0$ and
$D_\infty\ni\infty$. At the same time, it was well known that the
$\GRS$-method does not work in such a disconnected situation. In
the case of multipoint ($m$-point) PA the situation is even more
complicated, since it should be treated as a problem of an optimal
partition of the Riemann sphere into $m$ domains. Nevertheless,
the suitable generalization and improvement of the original
version of the $\GRS$-method to this situation was given by
Buslaev in~\cite{Bus13}. This generalization provided Buslaev the
opportunity to extend Stahl theory to multipoint PA for
multivalued analytic functions.\footnote{And even to extend the
Stahl theory to more general situation, namely, when the $m$-point PA
corresponds to a set of germs
$\myf=\{f_1,\dots,f_m\}=\{(f_1,z_1),\dots,(f_m,z_m)\}$ of $m$
multivalued analytic functions, each of which has a finite set of
singular points on the Riemann sphere. These $m$ germs are given
at $m$ distinct points $z_1,\dots,z_m$ and it may happen that not
even one of these germs can be obtained via the analytic
continuation of another germ.} In fact, Buslaev made a new step towards
the development of the $\GRS$-method in order to extend it to the
disconnected case. Ultimately, it made the $\GRS$-method much more
powerful than it was before.

\subsection{Buslaev's \texorpdfstring{$S$}{S}-curve
and Buslaev's limit zero-pole distribution theorem
in two-point Pad\'e approximants theory}\label{s3s1}
%\subsection{Two-point Pad\'e approximants: Buslaev's \texorpdfstring{$S$}{S}-curve
%and Buslaev's limit zero-pole distribution theorem}\label{s3s1}

Let us discuss in short the case of multipoint PA and Buslaev Theorem. The
original version of Buslaev Theorem is of general kind,
since it deals with $m$-point PA.
To be more precise, given a set of $m$ distinct points
$z_1,\dots,z_m\in\CC$ and a set $\myf=\{f_1,\dots,f_m\}$ of $m$ germs of $m$ analytic functions,
such that $f_j\in\HH(z_j)$, $j=1,\dots,m$, we seek
a rational function $B_n\in\CC_n(z)$ in $z$ of order $\leq{n}$, and such that
the following relations\footnote{The corresponding relation should be changed
when $z_j=\infty$ for some $j\in\{1,\dots,m\}$, see~\eqref{b5} below.}
hold:
\begin{equation}
(f_j-B_n)(z)=O\bigl((z-z_j)^{n_j}\bigr),\quad z\to z_j,
\label{b1}
\end{equation}
where $\sum\limits_{j=1}^m n_j=2n+1$, $n_j\in\ZZ_{+}$, $j=1,\dots,m$.
Suppose that all functions $f_j$ are multivalued analytic
functions with finite sets of branch points. Generically, all
functions $f_j$ are distinct, that is, not even one of $f_j$ is an analytic
continuation of another, say $f_k$, $k\neq j$.
% Suppose additionally,
% preserving the conditions of Buslaev Theorem, that $n_j/n\to2p_j$ as
% $n\to\infty$, where $\sum_{j=1}^mp_j=1$, $p_j\geq0$.

Suppose that each $f_j\in\myA^0(\myo\CC\setminus\Sigma_j)
:=\myA(\myo\CC\setminus\Sigma_j)\setminus\HH(\myo\CC\setminus\Sigma_j)$,
where $\Sigma_{j}=\Sigma_{f_j}$ and $\#\Sigma_j<\infty$, that is, each
$f_j$ is a multivalued analytic function on the Riemann sphere,
which is punctured at a finite set of branch points of~$f_j$. Suppose that
in~\eqref{b1}, preserving the conditions of Buslaev Theorem,
$n_j/n\to2p_j$ as $n\to\infty$, where
$\sum_{j=1}^mp_j=1$, $p_j\geq0$. Buslaev Theorem states that there exists, in the generic
case, an ``optimal partition'' of the Riemann sphere into $m$~domains
$D_j\ni z_j$, which are ``centered'' at the given interpolation points $z_j$, and
such that the compact set
$\Gamma:=\myo\CC\setminus\bigsqcup_{j=1}^mD_j$ consists of a finite
number of analytic curves and possesses the property of ``symmetry'' (see~\eqref{b3}).
The compact set $\Gamma$ itself is an $S$-curve ``weighted'' in the
external field, given by the logarithmic potential of the negative unit
charge $-\nu$, concentrated at the set of the given points
$\{z_1,\dots,z_m\}$:
\begin{equation}
\nu=\sum_{j=1}^mp_j\delta_{z_j}.
\label{bus3}
\end{equation}
Furthermore, each function $f_j\in\HH(z_j)$ possesses a holomorphic
(analytic and single-valued) continuation into $D_j$, $f_j\in\HH(D_j)$,
$j=1,\dots,m$. Finally, the $m$-point PA $B_n$ converges in
capacity, as $n\to\infty$, inside each domain $D_j$ to the function
$f_j$,
\begin{equation}
B_n(z)\overset{\mcap}\longrightarrow f_j(z),\quad n\to\infty,
\quad z\in D_j.
\label{b2}
\end{equation}
There exists a limit zero-pole distribution of the multipoint PA $B_n$.
The distribution in question coincides with the probability measure $\lambda=\lambda_\Gamma$,
which is concentrated on~$\Gamma$ and is equilibrium in the external
field $V^{-\nu}$. This field is determined by the negative charge concentrated on the set
$\{z_1,\dots,z_m\}$ of the interpolation points, each having the ``weight'' of $-p_j$ at
each point $z_j$ (see~\eqref{bus3}). Also, the $S$-property of the curve
$\Gamma$ and the rate of convergence in~\eqref{b2} may be
completely characterized by (cf.~\eqref{1.6} and \eqref{1.6eq}):
\begin{equation}
\frac{\partial (V^\lambda-V^\nu)}{\partial n^{+}}(z)=
\frac{\partial (V^\lambda-V^\nu)}{\partial n^{-}}(z),
\quad z\in\Gamma^0,
\label{b3}
\end{equation}
and (cf.~\eqref{2.3})
\begin{equation}
\bigl|f_j(z)-B_n(z)\bigr|^{1/n}\overset{\mcap}\longrightarrow
e^{-2\sum\limits_{k=1}^m p_k g^{}_{D_k}(z,z_k)},\quad n\to\infty,\quad z\in D_j,
\label{b4}
\end{equation}
where $g_{D_k}(z,z_k)$ is the Green's function for the domain $D_k$, with a
singularity at the point $z=z_k\in D_k$, $k=1,\dots,m$.

In what follows, for the sake of simplicity, we restrict our attention to the
particular case $m=2$ of Buslaev Theorem. Thus, we will discuss
in details only the case of two-point Pad\'e approximant.

Let $z_1=0$, $z_2=\infty$ and $\myf=\{f_0,f_\infty\}$ be the set of two
multivalued analytic functions, such that $f_0\in\HH(0)$ and
$f_\infty\in\HH(\infty)$, and also
$f_0,f_\infty\in\myA^0(\myo\CC\setminus\Sigma)$, where
$\#\Sigma<\infty$. Thus, each of the functions $f_0$ and $f_\infty$ is a
multivalued analytic function on the Riemann sphere, punctured at a
finite set of points, each point is a branch point of $f_0$ or
of $f_\infty$ or of both of them. In other words, $f_0$ and $f_\infty$ are
the two germs of the multivalued analytic functions, given at the points
$z_1=0$ and $z_2=\infty$, respectively. It is worth noting that they may be
two germs of the same analytic function, taken at two different points,
namely at $z=0$ and at $z=\infty$.

The two-point\footnote{In the classical terminology, this is the $n$-th truncated fraction
of the classical $T$-fraction, see~\cite{JoTh85}, and also~\cite{Bus01}.} PA is defined as follows.
Given $n\in\NN$, let $P_n,Q_n\in\CC_n(z)$, $Q_n\not\equiv0$, be polynomials
of degree $\leq{n}$, such that\footnote{For a fixed $n\in\NN$,
we can also claim that the left side of~\eqref{b5} is $O(z^{n+1})$
as $z\to0$ and $O(1)$ as $z\to\infty$, but this does not change the
convergence theorem itself.} the following relations hold
\begin{equation}
R_n(z):=\bigl(Q_nf-P_n\bigr)(z)=\begin{cases}
O(z^n),& z\to0,\\
O(1/z),& z\to\infty.
\end{cases}
\label{b5}
\end{equation}
The pair of polynomials $P_n$ and $Q_n$ is not unique, but the rational
function $B_n:=P_n/Q_n$ is uniquely determined by~\eqref{b5},
and is called the two-point diagonal PA to the set of germs of the
functions~$\myf=\{f_0,f_\infty\}$.
In the generic case, from~\eqref{b5} follows that
\begin{equation}
\bigl(f-B_n\bigr)(z)=\begin{cases} O(z^{n}),&z\to0,\\
O(1/z^{n+1}),& z\to\infty.
\end{cases}
\label{b6}
\end{equation}
If it exists, the rational function $B_n=B_n(z;f)\in\CC_n(z)$ is uniquely
determined by the relation~\eqref{b6}.

As for the classical Stahl's case, the existence of an $S$-curve,
associated with the two-point PA and weighted in the external field $V^{-\nu}$,
$\nu=(\delta_0+\delta_\infty)/2$, is the crucial element for
Buslaev's two-point convergence theorem. Such a weighted $S$-curve
$\Gamma=\Gamma_{\bus}(f_0,f_\infty)$ exists\footnote{In general,
there may exist some degenerated cases.} (see~\cite{BuMaSu12})
and makes an ``optimal'' partition of the Riemann
sphere into two domains $D_0\ni0$ and $D_\infty\ni\infty$, such that
$\myo\CC=D_0\sqcup\Gamma\sqcup D_\infty$, $f_0\in\HH(D_0)$ and
$f_\infty\in\HH(D_\infty)$. The compact set $\Gamma$ is a weighted $S$-curve, i.e.
$\Gamma$ consists of a finite number of analytic arcs and possesses the
following property of ``symmetry'':
\begin{equation}
\frac{\partial(V^\lambda-V^{\nu})}{\partial n^{+}}(z)
=\frac{\partial(V^\lambda-V^{\nu})}{\partial n^{-}}(z),
\quad z\in \Gamma^0,
\label{b7}
\end{equation}
where $\lambda=\lambda_\Gamma$ is the probability measure
concentrated on~$\Gamma$ and the equilibrium measure in the external field
$V^{-\nu}(z)=\frac12\log|z|$; furthermore, $\lambda$ is generated by the negative unit charge
$-\nu$, that is,
\begin{equation}
V^\lambda(z)-\frac12\log|z|\equiv\const,\quad z\in\Gamma.
\label{b8}
\end{equation}
As before, $\Gamma^0$ is the union of all open arcs of~$\Gamma$ (the
closures of which constitute $\Gamma$) and $\partial n^{+}$ and $\partial n^{-}$
are the inner (with respect to $D_0$ and $D_\infty$) normal
derivatives at a point $z\in \Gamma^0$ from the opposite sides of
$\Gamma^0$. Clearly, $\lambda$ is the balayage of the measure $\nu$ from
$D_0\sqcup D_\infty$ onto~$\Gamma$. It is worth noting that $\Gamma$
itself is a union of the closures of the critical trajectories of a
quadratic differential (see~\cite{BuMaSu12},
\cite{Bus13},~\cite{KoSu14}).

%\begin{remark}\label{rem3}
Here, for the sake of simplicity, we only consider the case of two-point PA,
and we set $z_1=0$ and $z_2=\infty$. In what follows, we also suppose that
$f_0$ and $f_\infty$ are the germs of the same multivalued analytic function
$f$, and denote them by $f_0\in\HH(0)$ and $f_\infty\in\HH(\infty)$.
We suppose that the function $f$ has a finite set of singular points in
$\myo\CC$. For example, $f$ may be an algebraic function, i.e. a function
given by an algebraic equation over the field $\CC(z)$, or a solution
of a linear homogeneous differential equation with polynomial coefficients
from the ring $\CC(z)$
(see~\cite{Lag85},~\cite{Per57},~\cite{ChCh84},~\cite{MaRaSu12}).

\begin{figure}[!ht]
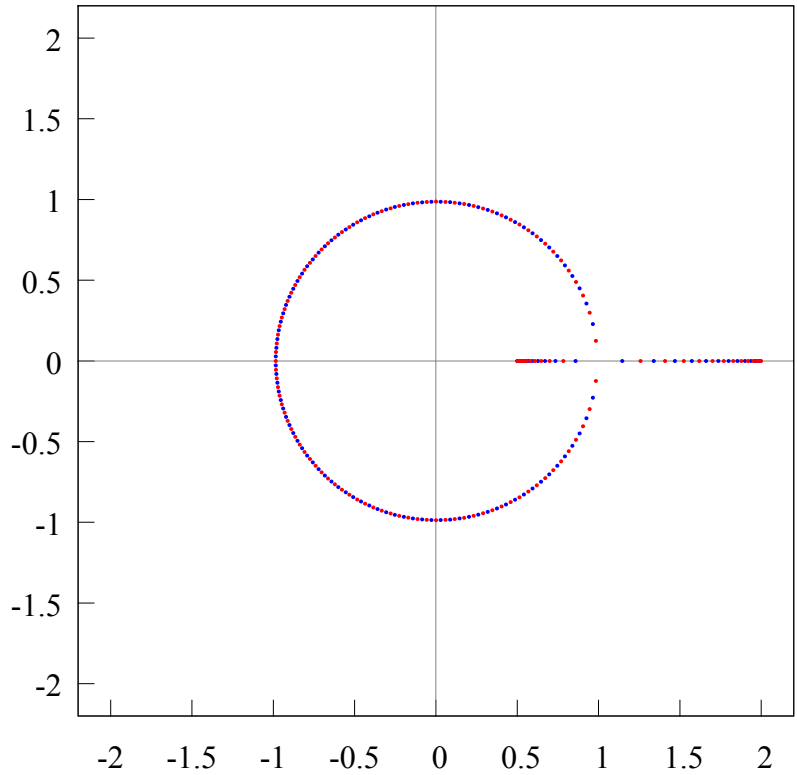

\centerline{
%\inclgr{Pic/buslaev_compact}}
\inclgr{FIG2_1}}
\vskip-6mm
\caption{Numerical zeros and poles distribution of two-point Pad\'e
approximants $[120/120]$ to the set of functions $\myf=\{f_0,f_\infty\}$,
where $f_0=((1-2z)(2-z))^{-1/2}$,
$f_\infty=((2z-1)(z-2))^{-1/2}+1$. The germs $f_0$ and $f_\infty$ result in
two different multivalued analytic functions. Thus, this is a generic case
and by Buslaev Theorem the associated $S$-curve partitions
the Riemann sphere into two domains.}
\label{Fig_bus210b_4000_120_full}
\end{figure}

Notice that the functions $f_0(z)=(1-z^2)^{-1/2}\sim1$, $z\to0$, and
$f_\infty=(z^2-1)^{-1/2}\sim1/z$, $z\to\infty$, are the germs of the
same analytic function~$f$, given by the equation $(z^2-1)w^2=1$. But the
functions $f_0(z)=(1-z^2)^{-1/2}$ and $f_\infty=(z^2-1)^{-1/2}+1$ are
not so. Thus, the latest case is the generic case, and hence $D_0\cap
D_\infty=\varnothing$ (see Fig.~\ref{Fig_bus210b_4000_120_full}).
%\end{remark}

Now we are ready to formulate the particular case of Buslaev
Theorem for two-point PA (cf. Stahl Theorem).

\begin{theoremBusTwoPoint}[V.~I.~Buslaev~\cite{Bus13}]
Let the function $f\in\HH(0)\cap\HH(\infty)$,
$f\in\myA^0(\myo\CC\setminus\Sigma)$, $\#\Sigma<\infty$, and let the pair
of germs $f_0,f_\infty$ be in a common
position\footnote{Equivalently, we say that Buslaev's
$S$-curve $\Gamma$ divides the Riemann sphere into two domains.}.
Let $\myo\CC=D_0\sqcup\Gamma\sqcup D_\infty$ be the optimal partition of
the Riemann sphere into two domains $D_0\ni0$ and $D_\infty\ni\infty$,
such that $f_0\in\HH(D_0)$, $f_\infty\in\HH(D_\infty)$,
$D_0\cap D_\infty=\varnothing$, and $\Gamma$
possesses the weighted $S$-property with respect to the external field
$V^{-\nu}$, $\nu=(\delta_0+\delta_\infty)/2$. Then for the $n$-diagonal
two-point PA $B_n$ of the set of the germs $\myf=\{f_0,f_\infty\}$ the following
statements hold true:

1) there exists a limit zero-pole distribution for $B_n$, namely,
\begin{equation}
\frac1n\chi(P_n),\frac1n\chi(Q_n)\overset{*}\longrightarrow\lambda_\Gamma,
\quad n\to\infty;
\label{b9}
\end{equation}

2) there is convergence in capacity as $n\to\infty$, namely,
\begin{equation}
B_n(z)\overset{\mcap}\longrightarrow f_0(z),\quad z\in D_0,\quad
B_n(z)\overset{\mcap}\longrightarrow f_\infty(z),\quad z\in D_\infty;
\label{b10}
\end{equation}

3) the rate of convergence in~\eqref{b10} as $n\to\infty$ is completely
characterized by the relations
\begin{equation}
\begin{aligned}
\bigl|f_0(z)-B_n(z)\bigr|^{1/n}
&\overset{\mcap}\longrightarrow
e^{-g^{}_{D_0}(z,0)},\quad z\in D_0,\\
\bigl|f_\infty(z)-B_n(z)\bigr|^{1/n}
&\overset{\mcap}\longrightarrow
e^{-g^{}_{D_\infty}(z,\infty)},\quad z\in D_\infty.
\end{aligned}
\label{b11}
\end{equation}
\end{theoremBusTwoPoint}

\begin{remark}\label{rem5}
As in Stahl Theorem, item~1) is equivalent to the following
relation as $n\to\infty$ for the monic two-point Pad\'e
polynomials:
\begin{equation}
|P^*_n(z)|^{1/n},|Q^*_n(z)|^{1/n}
\overset{\mcap}\longrightarrow e^{-V^\lambda(z)},\quad z\in D_0\sqcup
D_\infty.
\label{b12}
\end{equation}
\end{remark}

\begin{remark}\label{rem6}
``Optimal'' should be understood in connection with the
$S$-property~\eqref{b7}. In fact, this property means that the compact
set~$\Gamma$ possesses a ``stationary'' or ``equilibrium'' property in
the presence of the external field $V^{-\nu}$. This stationary property
is of unstable type, see~\cite{Rak12}.

In Buslaev Two-Point Theorem, equality~\eqref{b7} is more
complicated than equality~\eqref{1.6} in Stahl Theorem. To be more
precise, equality~\eqref{b7} should be understood as follows. In
the generic case, the compact set $\Gamma$ divides the complex
plane into two domains $D_0\ni0$ and $D_\infty\ni\infty$, such
that all three sets $\gamma:=\partial D_0\cap\partial D_\infty$,
$\gamma_0:=\partial D_0\setminus\partial D_\infty$ and
$\gamma_\infty:=\partial D_\infty\setminus\partial D_0$ are
nonempty sets. When $z\in\gamma$, $\partial n^{-}$ is the
normal derivative at the point $z$ from the boundary to the inside
of $D_0$ and $\partial n^{+}$ is the normal derivative at the
point $z$ from the boundary to the inside of $D_\infty$. If
$z\in\gamma_0$, $\partial n^{-}$ and $\partial n^{+}$ are the
normal derivatives at the point $z$ to the inside of $D_0$ from
the opposite sides of $\gamma_0$. We treat the case
$z\in\gamma_\infty$ in a similar way, but with respect to the
domain~$D_\infty$.

It should be emphasized that the main problem in this direction is to
prove the existence of a stationary compact set $\Gamma$, and
to characterize it as a weighted $S$-curve. The fact, that the equilibrium
measure $\lambda$ is the balayage of the measure $\nu$ from $D_0\cup
D_\infty$ onto $\Gamma$, is a trivial one.
\end{remark}

\begin{remark}\label{rem7}
We regard the ``optimal'' partition as optimal with respect to the given
field\footnote{Here the potential $V^{-\nu}(z)=\frac12\log|z|$ is spherically normalized.}
 $V^{-\nu}$, where
$\nu=\frac12(\delta_0+\delta_\infty)$. In general, in the two-point case
we have $\nu=p_0\delta_0+p_\infty\delta_\infty$, where $p_0,p_\infty>0$,
$p_0+p_\infty=1$. Thus, for the fixed function
$f\in\myA^0(\myo\CC\setminus\Sigma)$ and the fixed points $z_0=0$ and
$z_\infty=\infty$, the optimal partition and the compact set $\Gamma$
depend on the pair $p_0,p_\infty$
(see~\cite{Bus13}).
\end{remark}

\subsection{Strong asymptotics in two-point Pad\'e approximants theory}\label{s3s2}

It is worth noting that, in general, the $\GRS$-method is far from
being complete, because in the theory of LZD of HP polynomials
some complicated theoretical potential problems arise in a natural
way. At the moment these problems cannot be solved via Buslaev's
approach\footnote{In Buslaev Theorem, a set $\myf$ of germs of
multivalued functions $f_j$ should satisfy a more restricted
conditions, than in the original Stahl Theorem. Namely, $\myf$
should have a finite set of singular points instead of a set of
zero capacity, as it supposed in Stahl Theorem.}
(see~\cite{RaSu12},~\cite[Sec~6.2]{RaSu13}). The reason is that in
the case of the $m$-point PA, the external field is generated by
the finite number of $m$ positive pointed charges, concentrated at
the $m$ interpolation points (see
also~\cite{BaStYa12},~\cite{BuSu14b}). In contrast to this fact,
 in the case of
HP polynomials
the external field is of more complicated structure (see for example~\cite{RaSu13}).

In the case of multipoint PA (in part, in the case of two-point PA),
the situation is similar to the case of the
classical PA. Namely, given Buslaev's compact set (i.e. weighted
$S$-curve) $\Gamma_{\bus}$, the associated canonical two-sheeted
Riemann surface immediately appears, and it is equipped with the
corresponding Abelian integrals, etc. Recently, A.~V.~Komlov and
S.~P.~Suetin~\cite{KoSu14} derived a formula for the strong asymptotics
of the two-point Pad\'e polynomials for a special class of multivalued analytic
functions, given by the representation
\begin{equation}
f(z)=\(\frac{z-a_1}{z-a_2}\)^\alpha,
\quad \alpha\in\CC\setminus\QQ .
\label{2.8}
\end{equation}
The result was obtained by use of the following reasoning:
first, the authors proved that the corresponding two-point Pad\'e polynomials of
degree $n$ and the corresponding remainder function are just the
independent solutions of a linear differential equation of second order,
which contains some polynomial accessory parameter of fixed degree, but
depending on~$n$. Second, they proved that Buslaev's compact set,
associated with the function~\eqref{2.8}, yields the Stokes
lines\footnote{Since the accessory parameter depends on~$n$, it would be
better to say that Buslaev's compact set attracts
the Stokes lines of the differential equation, as $n\to\infty$.} for this
differential equation. Finally, following Nuttall's
approach~\cite{Nut86} (see also~\cite{MaRaSu12}), a suitable
modification of the classical Liouville--Steklov method
(in accordance with the presence of accessory parameter)
was used to derive formulae of strong asymptotics of
two-point Pad\'e polynomials\footnote{It is natural to call
these polynomials two-point Jacobi polynomials, see~\eqref{ort-two-point}.} and the remainder
function. It follows from the former and the latter, that the existence of Stahl's
$S$-curve, in the case of the classical PA, and the existence of
Buslaev's weighted $S$-curve, in the case of two-point PA, play the
crucial role not only in the problem of limit zero-pole distribution of
PA, but also in the problem of strong asymptotics of PA in both cases.
Thus, for the moment, the problem of weak-star asymptotics is completely
solved in both cases (for classical PA and for $m$-point PA it was completed by
H.~Stahl~\cite{Sta85a}--\cite{Sta86b} and by
 V.~I.~Buslaev~\cite{Bus13}, \cite{BuMaSu12},~\cite{Bus15}, respectively).
These results should be regarded as successful applications of the
$\GRS$-method. Notice that the first result is connected with the pure
logarithmic equilibrium problem, without any external field, and the
second one is connected with logarithmic equilibrium problem under the
presence of the external field. In general, the problem of
strong asymptotics of the Pad\'e polynomials is still open in both cases. However,
for the moment, some promising results are obtained in both cases, and all
of them are based on the existence of the associated $S$-curve, and
the associated weighted $S$-curve as well.

Notice that the two-point Pad\'e polynomial $Q_n$ of the
function~$f$, given by~\eqref{2.8}, satisfies the following
non-Hermitian orthogonality relations, with respect to the variable
weight function $(\zeta-a_1)^\alpha(\zeta-a_2)^{-\alpha}/\zeta^n$,
\begin{equation}
\label{ort-two-point}
\oint_\gamma \zeta^k \(\frac{\zeta-a_1}{z-a_2}\)^\alpha
\frac{d\zeta}{\zeta^n}=0,\qquad k=0,\dots,n-1,
\end{equation}
where $\gamma$ is an arbitrary contour, that creates a path through the points $a_1,a_2$,
and separates the point $z=0$ from the infinity point~$z=\infty$.

In the general theory of HP polynomials, the situation is quite different
from above. In that direction, the first results of general
character were obtained by A.~A.~Gonchar and E.~A.~Rakhmanov~\cite{GoRa81} in
1981 for the case of pure Angelesco systems of Markov-type functions. In
1984, J.~Nuttall~\cite{Nut84} stated some conjectures about HP polynomials,
that were important and with an impact to the theory of HP
polynomials. In 1986, E.~M.~Nikishin~\cite{Nik86} discovered and investigated
new systems of Markov-type functions, that are still of great
interest for all experts in HP polynomials theory.\footnote{One of the
main reasons is that in the generic case the pair of functions
$f,f^2$ form a Nikishin system.} These systems are named after him,
Nikishin systems. In 1988, H.~Stahl~\cite{Sta88} published a
survey,\footnote{The full manuscript of~\cite{Sta88} was available
much earlier than 1988. It is much bigger than the length of the published paper,
and is now available in electronic form.} where he submitted a lot of
conjectures about LZD of type~I and type~II HP polynomials, that
were based, in part, on the paper by Gonchar and Rakhmanov~\cite{GoRa81}.
As in Nuttall's paper~\cite{Nut84}, the conjectures from~\cite{Sta88}
were about the real and the complex case. In 1997, A.~A.~Gonchar,
E.~A.~Rakhmanov, and V.~N.~Sorokin~\cite{GoRaSo97} proved a general result about LZD of type~II
HP polynomials for the mixed, i.e. Angelesco and Nikishin, systems of
Markov-type functions. In 2010, A.~I.~Aptekarev and V.~G.~Lysov~\cite{ApLy10}
proved the most general at the moment result about LZD of type~II HP
polynomials for the mixed systems of Markov-type functions. To be more exact, they imposed
 the case when the support of one function of Markov type among the
system under consideration is a {\it proper part} of the support of another
Markov-type function among this system (for the
proof of a partial case via another approach the reader is referred to~\cite{Rak11}). All these general
results~\cite{GoRa81},~\cite{GoRaSo97},~\cite{ApLy10} are
dealing with systems of Markov-type functions with supports lying on the real
line. But not even one of them may be applied to the
situation under question in the current paper: when the
supports of two Markov-type functions have \textit{nonempty intersection, but not
even one of them is a proper part of the other}.

In the fundamental paper by Nuttall~\cite{Nut84} were proposed several conjectures of
a general type about asymptotics of HP polynomials.
All of them are concerned with both type~I and type~II HP polynomials,
and only with the strong asymptotics of HP
polynomials. Neither the existence of the associated $S$-curve, nor
the weak-star asymptotics of HP polynomials were discussed
in~\cite{Nut84}. This situation is typical for all results of HP polynomials. In
fact, the only results of general character on this subject were
obtained for the ``real case'' situation. In other words, all rigorous
and general results in this direction were proved for the
systems of Markov-type functions, and under such the conditions that LZD
of HP is completely described by the associated ``matrix of interaction''
and extremal theoretical potential problem. Depending on the system of
Markov-type functions, this matrix may be of Angelesco type, or of
Nikishin type, or some kind of ``mixed'' type. In each case,
the problem of LZD of HP polynomials is solved in terms of
the corresponding equilibrium measures, concentrated on {\it finite number
of segments of the real line}. The most general description of such ``real
case'' was done by Gonchar and Rakhmanov~\cite{GoRa85} (see also the
above cited papers~\cite{GoRa81},
\cite{Nik86},~\cite{GoRaSo97},~\cite{ApLy10}). It is worth noting
that some special properties of the HP polynomials in the real case
were proved by G.~L\'opez Lagomasino and
coauthors~\cite{FiLo11},~\cite{DeLoLo13},~\cite{LoMeFi15}, namely,
the property of normality, the convergence property (but without any
characterization of the rate of convergence), the interpolation
property, etc.

Thus, in all cases discussed above, the answer to the problem of LZD
of HP polynomials was given in terms in equilibrium measures
concentrated on the segments of the real line.\footnote{Because of this,
we refer to such type of case as ``real case''.} Presently, only a
few rigorous results are known of the asymptotics of HP polynomials for
the so-called ``complex case'' situation. In other words,
if the solution of the problem of LZD of HP polynomials is given in terms of
equilibrium measures concentrated on the associated $S$-curves, then the measures
are located somewhere in the complex plane, but not on the real line. As
usual, in such situations, we should first prove
the existence of an associated $S$-curve\footnote{In contrast
to Stahl's and Buslaev's theorems, in the case of HP for the collection
of three functions $[1,f_1,f_2]$, the associated
$S$-compact set should be considered as consisting of two
proper subsets, that play different roles in the associated theoretical
potential problem. The example of such $S$-compact set is given by
the ``Nuttall's condenser'',
see~\cite{RaSu12},~\cite{RaSu13},~\cite{KovSu14}.} (cf. Stahl's and
Buslaev's theorems). Also, only a few
rigorous results of LZD of HP polynomials are known
for the case when $m=3$, that is, for the
collection of three functions $[1,f_1,f_2]$ instead of the collection
of two functions $[1,f]$, as in the classical PA case. Notice that, for
the moment, there are two different approaches to the problem under
question. The first one is based on the cubic equation, and the other is
based on the concept of Nuttall's condenser. Here, we do not discuss the
details, but instead refer the reader to the
papers~\cite{NuTr87}, \cite{ApKuAs07},~\cite{ApTu12b},
\cite{ApTu14},~\cite{Apt08} and the references therein. For the cubic
equation and Nuttall's condenser, we refer
to~\cite{Kal79},~\cite{ApKa86},~\cite{ApKuAs07},
\cite{ApLyTu11},~\cite{MaRaSu15} and~\cite{RaSu12},~\cite{RaSu13},
\cite{Sue13b},~\cite{KovSu14}, respectively.

We can draw the following conclusion from the results of the papers reviewed above.
So far, there does not exist a general approach to the problem
of LZD of HP polynomials, and there is no connection to the powerful
$\GRS$-method, and, moreover, there does not even exist a suitable conjecture on this
subject, similar to Stahl's and Buslaev's results. In our knowledge,
all rigorously proved results are of partial character,
and there does not exist any theorem on LZD of HP polynomials,
with which we can explain from a theoretical
viewpoint the results of the numerical experiments, submitted in the
current paper. The numerical results, presented herein, can be described
in simple terms, but it is difficult to explain them from a
theoretical viewpoint. Ultimately, these results can be regarded as a
challenge to all experts on HP theory.

\section{Numerical results}\label{s4}

For completeness and for the reader's convenience, at first we present
some numerical results, concerned with Stahl's and Buslaev's theorems.

\subsection{Some numerical examples for classical Pad\'e approximants}\label{s4s1}

\subsubsection{The case of a function with three branch points,
Chebotarev--Stahl's $S$-curve and Stahl's limit zero-pole distribution
theorem}\label{s4s1s1}
\begin{figure}[!ht]
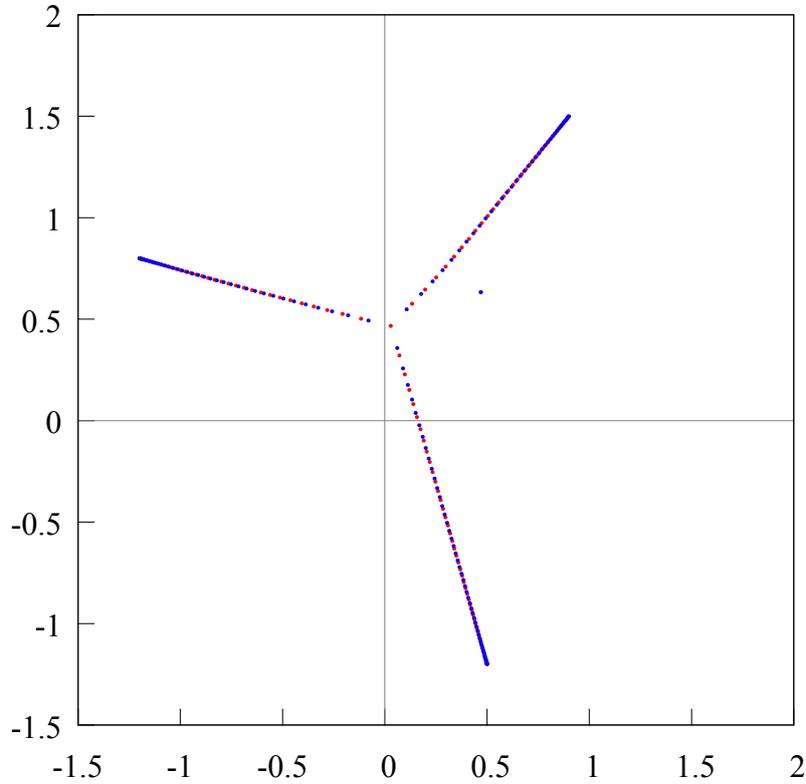

\centerline{
\inclgr{FIG1_1}}
\vskip-6mm
\caption{Zeros and poles of the diagonal Pad\'{e} approximant $[130/130]_f$ of the function
$f(z)=1/((z-(-1.2+0.8i))(z-(0.9+1.5i))
(z-(0.5-1.2i)))^{1/3}$, distributed accordingly
to the electrostatical model by E. A. Rakhmanov \cite{Rak12}.
There is a Froissart doublet (spurious zero-pole pair) when $n=130$ (see also
Fig.~\ref{Fig_pade10_2500_130_red} and Fig.~\ref{Fig_pade10_2500_130_blue}).
Since the genus of the Riemann surface is $1$,
there might be at most one Froissart doublet.
In full compliance with the Rakhmanov model \cite{Rak12},
the Froissart doublet ``attracts'' the Stahl $S$-compact $S_{130}$.
}
\label{Fig_pade10_2500_130}
\end{figure}

\begin{figure}[!ht]
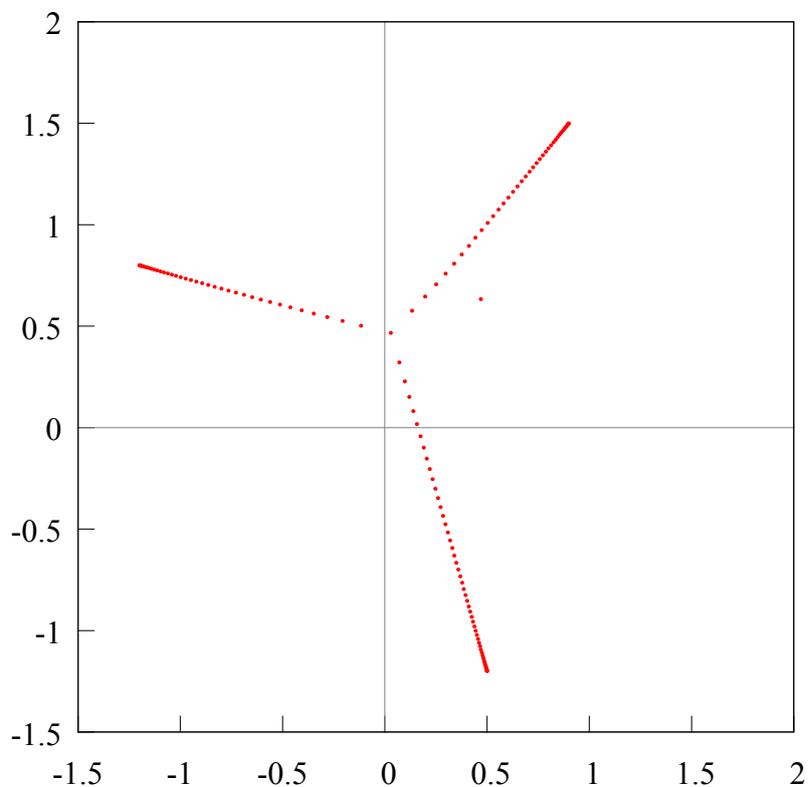

\centerline{
\inclgr{FIG1_2}}
\vskip-6mm
\caption{The poles of the Pad\'{e} approximant $[130/130]_f$
approximate a Chebotarev point $v_{130}$
for the $S$-compact $S_{130}$ (see \cite{Rak12}). The Chebotarev point is at $(0.029, 0.466)$.
When $n\to\infty$ we have that $v_n\to v$ is a classical
Chebotarev point. There is one spurious pole of the Pad\'{e} approximant $[130/130]_f$,
it is accompanied by a spurious zero of the Pad\'{e} approximant $[130/130]_f$
(see Fig.~\ref{Fig_pade10_2500_130_blue}).
The spurious zero-pole is at $(0.469, 0.633)$.
}
\label{Fig_pade10_2500_130_red}
\end{figure}

\begin{figure}[!ht]
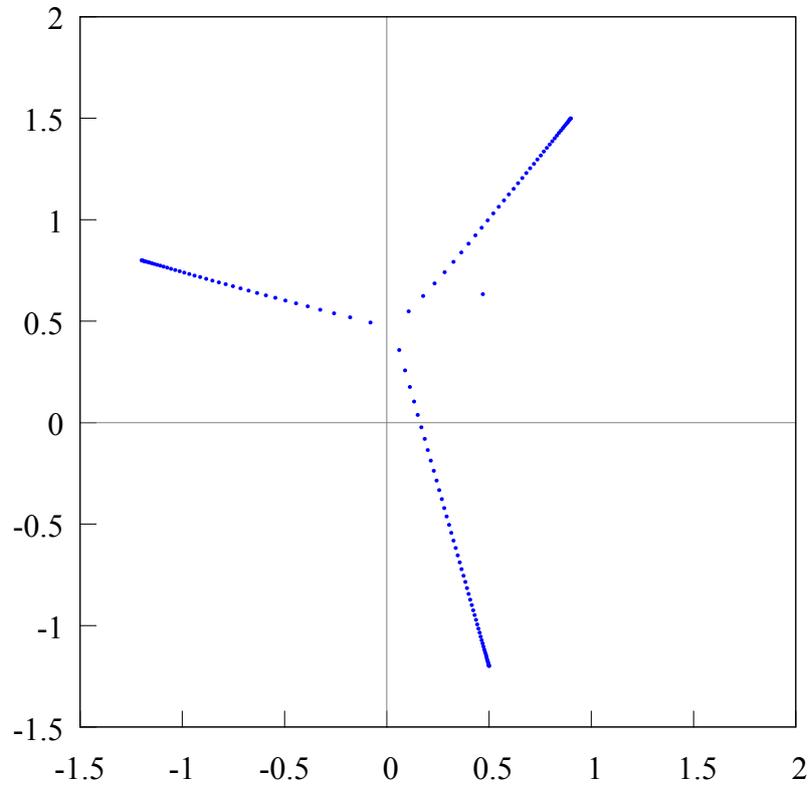

\centerline{
\inclgr{FIG1_3}}
\vskip-6mm
\caption{The Chebotarev point should not be approximated by
zeros of the Pad\'{e} approximant $[130/130]_f$ of the function \eqref{3.1}.
Evidently, the Cheboratev point does not exist on the picture.
There is one spurious zero of the Pad\'{e} approximant $[130/130]_f$,
it is accompanied by a spurious pole of the Pad\'{e} approximant $[130/130]_f$
(see Fig.~\ref{Fig_pade10_2500_130_red}).
}
\label{Fig_pade10_2500_130_blue}
\end{figure}

Let
\begin{equation}
f(z)=1/\bigl((z-a_1)(z-a_2)(z-a_3)\bigr)^{1/3},
\label{3.1}
\end{equation}
where $a_1=-1.2+0.8i$, $a_2=0.9+1.5i$, $a_3=0.5-1.2i$. On
Figure~\ref{Fig_pade10_2500_130}, the zeros (blue
points) and the poles (red points) of the PA $[130/130]_f$ of $f$ at infinity are plotted.
On Fig.~\ref{Fig_pade10_2500_130_blue} and
Fig.~\ref{Fig_pade10_2500_130_red} the zeros (blue points) and the
poles (red points) of PA $[130/130]_f$ are plotted separately.
Clearly, numerical zero-pole distribution is in good agreement
with the statements of Stahl Theorem. But there is a spurious
zero-pole pair, which does not correspond to any singularity of
the given function~\eqref{3.1}. The behavior of this pair as
$n\to\infty$ is not governed by Stahl Theorem, since this theorem is dealing with a
weak-star limit zero-pole distribution of PA. This kind of behavior was
discovered by Nuttall~\cite{Nut86} in 1986 via a suitable
modification of Liouville--Steklov method (see
also~\cite{MaRaSu12}). However, Nuttall's result was found only after Stahl Theorem was
proved. To be more precise, Nuttall used the existence of Stahl's
$S$-curve for the function~\eqref{3.1}.

Namely, Nuttall first
proved that Pad\'e polynomials of order~$n$ and the remainder
function for the function~\eqref{3.1} both solve a linear
homogeneous differential equation of second order, with polynomial
coefficients of fixed degrees. These polynomials have some
accessory parameters depending on~$n$. Subsequently, Nuttall
proved that Stahl's $S$-curve is the limit, as $n\to\infty$, of
the Stokes lines for these differential equations, that depend
on~$n$. Finally, in~\cite{Nut86} Nuttall proved that there may
exist only one spurious zero-pole pair of PA to the function~\eqref{3.1},
since the two-sheeted Stahl's Riemann surface $\RS_2$,
associated with the function~\eqref{3.1},
is of genus~$g=1$. Nuttall proved that if the zero and the pole of a
spurious zero-pole pair are close to one another, then they
actually cancel each other, as $n\to\infty$.
But if, for each $n\in\NN$, they do not coincide with
each other, then they are everywhere dense on the Riemann sphere,
as $n\to\infty$. Describing the behavior of the spurious zero-pole
pairs might require solving the so-called ``Jacobi inversion
problem'' (see~\cite{Spr57}), that is an equation with an Abelian
integral of first kind on the left and some expression that
is linear in~$n$ on the right. Thus, as $n\to\infty$, such
pairs form some type of ``winding of the torus'', which is
everywhere dense on that torus, that is, dense on the Stahl's
two-sheeted Riemann surface $\RS_2$ of genus~$g=1$. The fact that
the zero and the pole in some sense cancel each other in such a
zero-pole pair, as $n\to\infty$, is in full agreement
with the pure numerical results that were obtained by
M.~Froissart~\cite{Fro69}, see also~\cite[Chapter~2,
\S\,2.2]{BaGr81} and~\cite{Sue04}.

\subsubsection{The case of a function with several branch points,
Chebotarev--Stahl's $S$-curve and Stahl limit zero-pole distribution
theorem}\label{s4s1s2}
\begin{figure}[!ht]
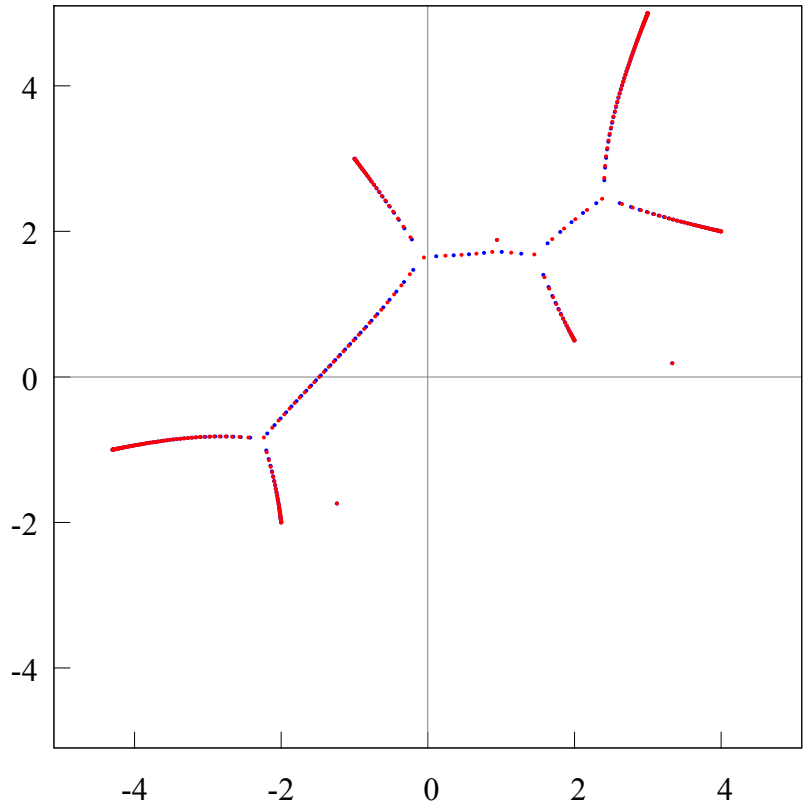

\centerline{
\inclgr{FIG1_4}}
\vskip-6mm
\caption{Zeros and poles of the diagonal Pad\'{e} approximant $[103/103]_f$ of the function
$f(z)=1/((z+(4.3+1.0i))
(z-(2.0+0.5i))(z+(2.0+2.0i))(z+(1.0-3.0i))
(z-(4.0+2.0i))(z-(3.0+5.0i)))^{1/6}$.
These zeros and poles are distributed in a plane, under fixed $n=103$,
accordingly to the electrostatical model by Rakhmanov~\cite{Rak12}.
Since the genus of the Riemann surface is $4$,
for each $n$ there might be no more than 4 Froissart doublets.
Here are observed 3 Froissart doublets.
In full compliance with the Rakhmanov model,
the Froissart doublets ``attract'' the Stahl $S$-compact $S_{103}$.
In general, the zeros and poles of the diagonal Pad\'{e} approximants $[n/n]_f$
are distributed as $n\to\infty$ accordingly to Stahl Theorem~\cite{Sta97b}.
}
\label{Fig_pade103_5000_266_full}
\end{figure}

Let
\begin{equation}
f(z)=1/\bigl((z-a_1)\cdot\dotso\cdot(z-a_6)\bigr)^{1/6},
\label{3.1.1}
\end{equation}
where $a_1=4.3+i$, $a_2=2+0.5i$, $a_3=2+2i$, $a_4=1-3i$,
$a_5=4+2i$, $a_6=3+5i$. On Fig.~\ref{Fig_pade103_5000_266_full}
are plotted zeros (blue points) and poles (red points) of the PA
$[103/103]_f$, taken at the infinity point $z=\infty$ to the
function~$f$ given by~\eqref{3.1.1}. Clearly, numerical zero-pole
distribution is in good agreement with the statements of Stahl
Theorem. But there are spurious zero-pole pairs, which do not
correspond to singularities of the given function~\eqref{3.1.1}.
The behavior of this pair, as $n\to\infty$, is not governed by
Stahl Theorem, since it is about weak-star limit zero-pole
distribution of PA. This behavior was discovered by
A.~Martinez-Finkelshtein, E.~A.~Rakhmanov, and
S.~P.~Suetin~\cite{MaRaSu12} in 2012 via a suitable modification
of Liouville--Steklov method (cf.~\cite{Nut86}). This result was
obtained on the base of Stahl Theorem. More precisely: the
authors of~\cite{MaRaSu12} used the basic fact about the existence
of Stahl's $S$-curve for the function~\eqref{3.1.1}. Namely, first
they proved that Pad\'e polynomials of order~$n$, as well as the
remainder function for the function~\eqref{3.1.1} solve a
linear homogeneous differential equation of second order with
polynomial coefficients of fixed degrees. These polynomials have
some accessory parameters depending on~$n$. Second, they
showed that Stahl's $S$-curve is the limit, as $n\to\infty$, of
the Stokes lines for these differential equations, depending
on~$n$. Ultimately, the authors concluded that
no more than four spurious zero-pole pairs of PA of the
function~\eqref{3.1.1} may exist.
The reason is that the Stahl's two-sheeted Riemann surface $\RS_2$,
which corresponds to the function~\eqref{3.1.1}, is of genus~$g=4$.
It was proved that the zero and the pole in each spurious zero-pole pair
are close to each other and they cancel each other, as $n\to\infty$.
However, in the ``generic case'', if they do not coincide with
each other for each $n\in\NN$, then they are everywhere dense on the Riemann sphere,
as $n\to\infty$. Describing the behavior of the spurious zero-pole
pairs might require solving a system of equations with some Abelian integrals
of first kind on the left side and some expressions, which are linear in~$n$, on the right side.
The fact that the zero and pole in some sense cancel each other in such a zero-pole pair
as $n\to\infty$ is in full agreement
with the pure numerical results obtained by M.~Froissart~\cite{Fro69},
see also~\cite[Chapter~2, \S\,2.2]{BaGr81} and~\cite{Sue04}.

\subsection{Some numerical examples for two-point Pad\'e approximants}\label{s4s2}

Let
\begin{equation}
f(z)=\(\frac{z-a_1}{z-a_2}\)^{1/4},
\label{3.2}
\end{equation}
where $a_1=0.9-1.1i$, $a_2=0.1+0.2i$.
Suppose that two ``different'' germs are taken, $f_0$ and $f_\infty$ of the
function~\eqref{3.2} at the point $z=0$ and at the infinity point
$z=\infty$, respectively. Here, ``different'' means that to
obtain the germ $f_\infty\in\HH(\infty)$ by the analytic continuation
of the germ $f_0\in\HH(0)$, we should go along a path that
encircles exactly one time one of the branch points ($a_1$ or $a_2$).

\begin{figure}[!ht]
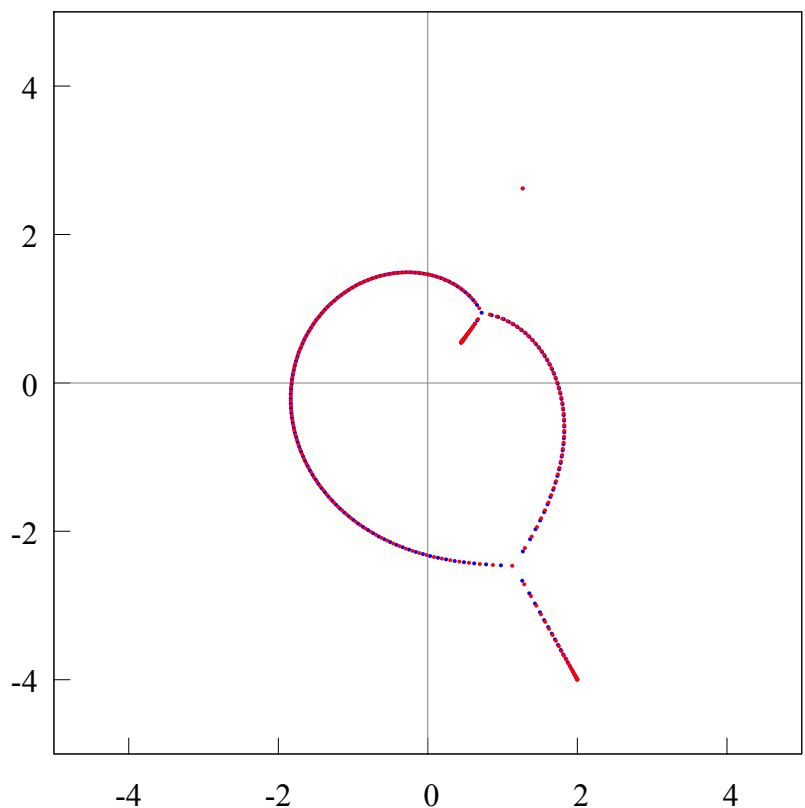

\centerline{
%\inclgr{Pic/buslaev_compact}}
\inclgr{FIG2_2}}
\vskip-6mm
\caption{Numerical zeros and poles distribution of two-point Pad\'e
approximants $[199/199]$ to the function $f(z)=((z-a_1)/(z-a_2))^{1/4}$,
where $a_1=0.9-1.1i$ and $a_2=0.1+0.2i$.
Here are selected two ``different branches'' of the function~$f$, namely,
$f_0=((z-a_1)/(z-a_2))^{1/4}$ and $f_\infty=-((z-a_1)/(z-a_2))^{1/4}$.
Almost all zeros (blue points) and poles (red points) approximate numerically
Buslaev's compact set. But there is a spurious zero-pole pair that moves
as $n\to\infty$ under the order of equation from~\cite{KoSu14}.}
\label{Fig_bus205c(2000)199_full}
\end{figure}

On Fig.~\ref{Fig_bus205c(2000)199_full} are plotted the zeros and poles of two-point PA
$[199/199]_{f_0,f_\infty}$ to the function~\eqref{3.2}. On Fig.~\ref{Fig_bus205c(2000)199_blue} and
Fig.~\ref{Fig_bus205c(2000)199_red} are plotted separately the zeros (blue points) and poles (red
points) of the two-point PA $[199/199]_{f_0,f_\infty}$ for the function~\eqref{3.2}.
This numerical zero-pole distribution of two-point PA
$[199/199]_{f_0,f_\infty}$ is in full agreement with Buslaev's limit
zero-pole distribution theorem~\cite{Bus13}. Similarly to the
function~\eqref{3.1} and classical PA, the two-sheeted Riemann surface,
associated with the function~\eqref{3.2} and the two-point PA, is of
genus~$g=1$. Therefore, similarly to the case of classical PA, there is a single
zero-pole pair of spurious character. The behavior of this pair, as
$n\to\infty$, cannot be described by Buslaev Theorem, since the theorem is only about
weak-star convergence of two-point PA.

The description of this behavior for the
special function~\eqref{3.2} and under the condition of a ``generic case'', imposed
on the branch points $a_1$ and $a_2$, was done
by Komlov and Suetin in~\cite{KoSu14}. Under this condition, in~\cite{KoSu14}
was derived a formula of the strong asymptotics for two-point PA.
The method producing this formula is similar to Nuttall's method for 
classical PA, see~\cite{Nut86} and~\cite{MaRaSu12}.
In~\cite{KoSu14}, it was proved first that for each
$n\in\NN$ the two-point Pad\'e polynomials and the corresponding
remainder function solve the homogeneous linear differential equation
of second order with polynomial coefficients of fixed degrees, but
depending on~$n$. After that, it was proved that Buslaev's $S$-curve
forms the limit, as $n\to\infty$, of the Stokes lines for these differential
equations. Finally, a formula of the strong asymptotics for the two-point PA
was found via the asymptotic Liouville--Steklov method.

\subsection{Some numerical examples for Hermite-Pad\'{e} polynomials}\label{s4s3}

Now we are ready to make some short description of the obtained
numerical results for the HP of the collection of three functions $[1,f_1,f_2]$.

\subsubsection{Parameter $a=-0.1$: logarithmic functions,
Fig.~\ref{Fig_log(pL1s)1500(200)_full}--Fig.~\ref{Fig_log(pL1s)1500(200)_blue},
square root functions,
Fig.~\ref{Fig_sqrt(pL1s)2000(200)_full}--Fig.~\ref{Fig_sqrt(pL1s)2000(200)_blue}
}

At first we set $a=-0.1$ (and keep $n=200$). Clearly, that value
of~$a$ is out of the range of $(0,1)$. But in that case all three pairs
$f_1,f_2$, $g_1,g_2$ and $h_1,h_2$ form
Angelesco systems, since $E_1,E_2\subset\RR$ and $E_1\cap
E_2=\varnothing$. On
Fig.~\ref{Fig_log(pL1s)1500(200)_full},
Fig.~\ref{Fig_log(pL1s)1500(200)_blue},
Fig.~\ref{Fig_sqrt(pL1s)2000(200)_full}, and
Fig.~\ref{Fig_sqrt(pL1s)2000(200)_blue}
are plotted the zeros of HP polynomials
$Q_{200,j}$ and $P_{200,j}$, $j=0,1,2$. The numerical
distribution of the zeros of the HP polynomials are similar in
both cases, and it does not depend on the different type of branching
of the functions $f_1,f_2$ and $g_1,g_2$ at the points $z=\pm1,\pm a$. The
reason is as follows. All these numerical results are in a full
agreement with the general theory of limit zero distribution (LZD) for
HP polynomials in the Angelesco case (see~\cite{GoRa81},
\cite{GoRaSo97},~\cite{ApLy10}; in these papers the case of type~II HP
polynomials was considered, but it is well known that in the Angelesco
case there is a direct connection between LZD of type~I and type~II HP
polynomials). Since the size of the supports of Markov-type functions $f_1$ and
$f_2$ (and $g_1$, $g_2$ as well) are equal, the phenomena of the
so-called ``pushing of the charge''\footnote{The phenomena was discovered
by Gonchar and Rakhmanov in~\cite{GoRa81}, 1981.} is absent in the
case when $a=-0.1$. It is also in full agreement with the
general theory (see~\cite{Kal79},~\cite{ApKa86},~\cite{GoRa81}). From
Fig.~\ref{Fig_log(pL1s)1500(200)_full} and
Fig.~\ref{Fig_log(pL1s)1500(200)_blue} follows, that the zeros of
HP polynomials $Q_{200,1}$ and $Q_{200,2}$ (and $P_{200,1}$,
$P_{200,2}$ as well) are located on the segments $E_1$ and $E_2$,
respectively. But the zeros of $Q_{200,0}$ (and $P_{200,0}$
respectively) are located on the imaginary axis. Thus, in some sense, the
zeros of $Q_{200,0}$ ``intend to separate'' the zeros of $Q_{200,1}$ and
$Q_{200,2}$ from each other. Again, this result is in good
agreement with the well known LZD of type~I HP polynomials in the
Angelesco case (see, e.g.~\cite{NuTr87}). To finish the analysis of the
case $a=-0.1$, we emphasize that in this case all numerical
results are in good agreement with the general theory of HP
polynomials. Thus it is reasonable to consider all numerical
results for $a=0.2;0.4;0.625;0.73;0.8$ as to be trustable.

\subsubsection{Parameter $a=0.2$: logarithmic functions,
Fig.~\ref{Fig_log(p2s)2000(200)_full}--Fig.~\ref{Fig_log(p2s)2000(200)_bk}
% square root functions,
% Fig.~\ref{Fig_sqrt(p2s)2000(200)_full}--Fig.~\ref{Fig_sqrt(p2s)2000(200)_bk},
% cubic root functions,
% Fig.~\ref{Fig_cub(p2s)2000(200)_full}--Fig.~\ref{Fig_cub(p2s)2000(200)_bk}
}

Now let $a=0.2$. For this value of the parameter $a$,
the intersection of the two segments $E_1=[-1,a]$ and
$E_2=[-a,1]$ equals $E_1\cap E_2=\Delta=[-a,a]\not=\varnothing$, that is, $\Delta$
is small, but nonempty real segment.

In all three cases, Case~\ref{cas1}, Case~\ref{cas2}, and
Case~\ref{cas3}, the numerical zero distribution (NZD) of HP
$Q_{200,j}(z;f_1,f_2)$, $P_{200,j}(z;g_1,g_2)$, and $\UU_{200,j}(z;h_1,h_2)$
for $j=1,2$ are similar to each other. Thus, here we only discuss the NZD for
HP $Q_{200,1}$ and $Q_{200,2}$ (these zeros are plotted by red and black
points respectively, see Fig.~\ref{Fig_log(p2s)2000(200)_red}
and~Fig.~\ref{Fig_log(p2s)2000(200)_bk}). Notice that zeros of $Q_{200,1}$ and $Q_{200,2}$ are
symmetric to each other, with respect to the imaginary axis. Thus, it is enough
to only consider the NZD for HP $Q_{200,2}$ (these zeros are plotted in black;
see~Fig.~\ref{Fig_log(p2s)2000(200)_bk}).

For this value of the parameter $a$, the essential part of zeros of HP
$Q_{200,2}$ are located not on the real line, but somewhere in the complex
plane. Evidently, this set is symmetric with respect to the real axis,
since we consider Markov-type functions. Due to the general
conjectures (see~\cite{Nut84},~\cite{Sta88}), that are based on the seminal
paper by Gonchar and Rakhmanov~\cite{GoRa81} about LZD of pure Angelesco
system of Markov-type functions, the LZD for HP $Q_{n,2}$, as $n\to\infty$,
should be described via an equilibrium measure, say $\lambda_2$, that
is associated with some special $\max$-$\min$ extremal theoretical
potential problem. From Fig.~\ref{Fig_log(p2s)2000(200)_bk} (see also
Fig.~\ref{Fig_log(p2s)2000(200)_red}) it follows immediately, that in the case
$a=0.2$ the support of the measure~$\lambda_2$ should be a disconnected
set. The same is true for the equilibrium measure $\lambda_1$, that
is associated with LZD for HP $Q_{n,1}$. The union $S_1(f_1,f_2)\cup S_2(f_1,f_2)$ of two
compact sets $S_1(f_1,f_2)=\supp\lambda_1$ and $S_2(f_1,f_2)=\supp\lambda_2$ form some
type of lenses. From the general approach, based on the paper by
Gonchar and Rakhmanov~\cite{GoRa87} about ``$1/9$-conjecture'', it follows, that the compact sets
$S_1(f_1,f_2)$ and $S_2(f_1,f_2)$ should be a weighted $S$-curve (see
also~\cite{Rak12}). Notice, that in both cases the open sets
$\myo\CC\setminus S_1(f_1,f_2)$ and $\myo\CC\setminus S_2(f_1,f_2)$ are also
domains (see Fig.~\ref{Fig_log(p2s)2000(200)_red} and
Fig.~\ref{Fig_log(p2s)2000(200)_bk} respectively).

It follows from Fig.~\ref{Fig_log(p2s)2000(200)_blue}, that the LZD for
HP $Q_{n,0}$ should be quite different from the LZD for HP $Q_{n,1}$
and $Q_{n,2}$. The reason is as follows. This LZD should be
described by an extremal $\max$-$\min$ theoretical potential problem of different type than
before. In part, the support $S_0(f_1,f_2)$ of the equilibrium
measure~$\lambda_0$, associated with this extremal theoretical potential problem, should be
a continuum, i.e. a connected compact set. But now, the open set
$\myo\CC\setminus S_0(f_1,f_2)$ is not a domain, since it consists of three domains.
The compact set $S_0(f_1,f_2)$ should also be a weighted $S$-curve, but of
some other nature than the compact sets $S_1(f_1,f_2)$ and $S_2(f_1,f_2)$. In particular,
there should be three Chebotarev's points on the compact set $S_0(f_1,f_2)$, i.e.
the points of zero density of the equilibrium measure~$\lambda_0$. One
of these points is located on the upper half plane, the other point
is located on lower half plane, and the third Chebotarev's point coincides
with the infinity point. Notice, that similarly to the case when
$a=-0.1$, now the $S$-curve $S_0(f_1,f_2)$ separates the $S$-curves $S_1(f_1,f_2)$ and
$S_2(f_1,f_2)$ from each other. It might be conjectured that the pair
$S_(f_1,f_2),S_2(f_1,f_2)$ forms some kind of weighted Nuttall's condenser
(see~\cite{RaSu13}).

\subsubsection{Parameter $a=0.2$: square root functions,
Fig.~\ref{Fig_sqrt(p2s)2000(200)_full}--Fig.~\ref{Fig_sqrt(p2s)2000(200)_bk},
and cubic root functions,
Fig.~\ref{Fig_cub(p2s)2000(200)_full}--Fig.~\ref{Fig_cub(p2s)2000(200)_bk}
}

The NZD of HP $P_{200,1}$ and $P_{200,2}$ for the square root functions,
given by~\eqref{1.3},~\eqref{1.4}, and NZD of HP polynomials $\UU_{200,1}$ and $\UU_{200,2}$
for the cubic root functions, given by~\eqref{c3.1},~\eqref{c3.2}, are similar to the NZD
of HP $Q_{200,1}$ and $Q_{200,2}$, for the pair of logarithmic
functions given by~\eqref{f1},~\eqref{f2}, respectively (see
Fig.~\ref{Fig_sqrt(p2s)2000(200)_red}, Fig.~\ref{Fig_sqrt(p2s)2000(200)_bk},
Fig.~\ref{Fig_cub(p2s)2000(200)_red}, and Fig.~\ref{Fig_cub(p2s)2000(200)_bk}).
It is also valid for the HP $\UU_{200,0}$,
i.e. this NZD is similar to the NZD of HP $Q_{200,0}$ for the given
logarithmic functions. Therefore, we do not discuss specially the NZD of HP
$P_{200,1}$, $P_{200,2}$, $\UU_{200,1}$, $\UU_{200,2}$, and $\UU_{200,0}$.
Notice, that there is a pair of spurious zeros of the HP $\UU_{200,0}$
(see Fig.~\ref{Fig_cub(p2s)2000(200)_blue}). From
that numerical fact follows, that the genus of the associated three-sheeted
Riemann surface $\RS_3$ should be equal to $2$ (cf.~\cite{NuTr87},~\cite{ApKuAs07}).

\subsubsection{Parameter $a=0.2$: square root functions, NZD of $P_{200,0}$,
Fig.~\ref{Fig_sqrt(p2s)2000(200)_blue}}

The NZD of HP polynomials $P_{200,0}$ is quite different from NZD of HP $Q_{200,0}$
(see Fig.~\ref{Fig_log(p2s)2000(200)_blue})
and $\UU_{200,0}$
(see Fig.~\ref{Fig_cub(p2s)2000(200)_blue}).
In fact, from the numerical results follows, that
the associated $S$-curve $S_0(g_1,g_2)$ should consist of two segments, and the
open set $\myo\CC\setminus S_0(g_1,g_2)$ should be a domain. In that case, the compact
set $S_0(g_1,g_2)$ separates distinctly the other $S$-curves $S_1(g_1,g_2)$
and $S_2(g_1,g_2)$ from each other.

\subsubsection{Parameter $a=0.4$:
logarithmic functions,
Fig.~\ref{Fig_log(p4s)2000(200)_full}--Fig.~\ref{Fig_log(p4s)2000(200)_bk},
square root functions,
Fig.~\ref{Fig_sqrt(p4s)2000(200)_full}--Fig.~\ref{Fig_sqrt(p4s)2000(200)_bk},
cubic root functions,
Fig.~\ref{Fig_cub(p4s)2000(200)_full}--Fig.~\ref{Fig_cub(p4s)2000(200)_bk};
parameter $a=0.625$:
logarithmic functions,
Fig.~\ref{Fig_log(p625s)2000(200)_full}--Fig.~\ref{Fig_log(p625s)2000(200)_bk},
square root functions,
Fig.~\ref{Fig_sqrt(p625s)2000(200)_full}--Fig.~\ref{Fig_sqrt(p625s)2000(200)_bk},
cubic root functions,
Fig.~\ref{Fig_cub(p625s)2000(200)_full}--Fig.~\ref{Fig_cub(p625s)2000(200)_bk}
}{\ }

For these two values of the parameter $a$ all the NZD of HP
$Q_{200,j}$, $P_{200,j}$, $\UU_{200,j}$, $j=0,1,2$, are similar to
the case when $a=0.2$. The observed difference does not change the
principal structure of the associated $S$-curves and is the following.
The lenses, that appeared for parameter $a=0.2$, become larger and larger when $a$
changes from $0.2$ to the values $a=0.4$ and $a=0.625$. The positions
of all four vertices of the lenses are not fixed, but depend on the value
of the parameter $a$. When the parameter $a$ increases from $a=0.2$ to
$a=0.4$ and after that to $a=0.625$, the two real vertices move from
the inside of the segment $[-1,1]$ towards to the end points $\pm1$. In
addition, the two pure imaginary vertices move along simultaneously from
the imaginary axis to the infinity point, where they meet each other,
under some critical value of the parameter $a^*\in(0.625,0.73)$ (see
the next Fig.~\ref{Fig_log(p73s)2000(200)_full}--Fig.~\ref{Fig_cub(p8s)2000(200)_bk};
the statement on the existence of the critical value $a^*$ is based only on
 numerical results and should be considered as a conjecture).
As usually, the NZD of HP $P_{200,0}$ corresponds to the whole segment $[-a,a]$.
Hence, the blue segment on
Fig.~\ref{Fig_sqrt(p2s)2000(200)_blue},
Fig.~\ref{Fig_sqrt(p4s)2000(200)_blue}, and
Fig.~\ref{Fig_sqrt(p625s)2000(200)_blue} becomes wider and wider as the
parameter $a$ increases form $a=0.2$ to $a=0.625$.

\subsubsection{Parameter $a=0.73$:
logarithmic functions,
Fig.~\ref{Fig_log(p73s)2000(200)_full}--Fig.~\ref{Fig_log(p73s)2000(200)_bk},
square root functions,
Fig.~\ref{Fig_sqrt(p73s)2000(200)_full}--Fig.~\ref{Fig_sqrt(p73s)2000(200)_bk},
cubic root functions,
Fig.~\ref{Fig_cub(p73s)1500(300)_full}--Fig.~\ref{Fig_cub(p73s)1500(300)_bk};
parameter $a=0.8$:
logarithmic functions,
Fig.~\ref{Fig_log(p8s)2000(200)_full}--Fig.~\ref{Fig_log(p8s)2000(200)_bk},
square root functions,
Fig.~\ref{Fig_sqrt(p8s)2000(200)_full}--Fig.~\ref{Fig_sqrt(p8s)2000(200)_bk},
cubic root functions,
Fig.~\ref{Fig_cub(p8s)2000(200)_full}--Fig.~\ref{Fig_cub(p8s)2000(200)_bk}
}

After $a>a^{*}$, where $a^{*}$ is the critical value of the parameter
$a$ described above, the NZD of HP $Q_{200,j}$, $P_{200,j}$,
$\UU_{200,j}$, $j=0,1,2$, dramatically changes. Namely, the complements
$\myo\CC\setminus S_1$ and $\myo\CC\setminus S_2$ of the $S$-curves
$S_1$ and $S_2$ are now disconnected open sets, instead of
domains, as they were when $a<a^{*}$. Plotted on a single
picture, all three numerical sets $S_1,S_2,S_3$, that consist of red,
black, and blue points, respectively, form in these three cases three quite different structures.
For instance, when the parameter $a$ equals $0.8$, in Case~\ref{cas1} the associated $S$-curve
$S_0(f_1,f_2)$ becomes a disconnected compact set, namely,
$S_0(f_1,f_2)$ consists of two continua (see
Fig.~\ref{Fig_log(p8s)2000(200)_blue}). In Case~\ref{cas2}, we have
$S_0(g_1,g_2)=[-a,a]$ (see Fig.~\ref{Fig_sqrt(p8s)2000(200)_blue}), and
finally, in Case~\ref{cas3} the $S$-curve $S_0(h_1,h_2)$ is a continua that
contains the infinity point (see Fig.~\ref{Fig_cub(p8s)2000(200)_blue}).
In the last case, $S_0(h_1,h_2)$ contains two Chebotarev's points
that are the points of zero density for the corresponding
equilibrium measure $\lambda_0$ with $\supp\lambda_0=S_0(h_1,h_2)$.

%%%endfulltext

\newpage
\clearpage

% LaTeX can keep in memory 18 floats at most
% \clearpage must be invoked each 18 graphics

%%%figures1

% Figures 1 - 4 (4 total)

%%%endfigures1

%%%figures2

% Figures 5 - 8 (4 total)

\begin{figure}[!ht]
\centerline{
%\inclgr{Pic/buslaev_compact}}
\inclgr{FIG2_3}}
\vskip-6mm
\caption{Numerical zeros distribution of two-point Pad\'e
approximants $[199/199]$ to the function $f(z)=((z-a_1)/(z-a_2))^{1/4}$,
where $a_1=0.9-1.1i$ and $a_2=0.1+0.2i$.
Here are selected two ``different branches'' of the function~$f$, namely,
$f_0=((z-a_1)/(z-a_2))^{1/4}$ and $f_\infty=-((z-a_1)/(z-a_2))^{1/4}$.
Almost all zeros (blue points) approximate numerically
Buslaev's compact set. But there is a spurious zero-pole pair that moves
as $n\to\infty$ under the order of equation from~\cite{KoSu14}.}
\label{Fig_bus205c(2000)199_blue}
\end{figure}

\begin{figure}[!ht]
\centerline{
%\inclgr{Pic/buslaev_compact}}
\inclgr{FIG2_4}}
\vskip-6mm
\caption{Numerical poles distribution of two-point Pad\'e
approximants $[199/199]$ to the function $f(z)=((z-a_1)/(z-a_2))^{1/4}$,
where $a_1=0.9-1.1i$ and $a_2=0.1+0.2i$.
Here are selected two ``different branches'' of the function~$f$, namely,
$f_0=((z-a_1)/(z-a_2))^{1/4}$ and $f_\infty=-((z-a_1)/(z-a_2))^{1/4}$.
Almost all poles (red points) approximate numerically
Buslaev's compact set. But there is a spurious zero-pole pair that moves
as $n\to\infty$ under the order of equation from~\cite{KoSu14}.}
\label{Fig_bus205c(2000)199_red}
\end{figure}

%%%endfigures2

%%%figures3

% Figures 9 - 12 (4 total)

\begin{figure}[!ht]
\centerline{
\inclgraph{FIG3_1}}
\vskip-6mm
\caption{Numerical distribution of zeros of type~I Hermite--Pad\'e polynomials
$Q_{200,0}$ (blue points),
$Q_{200,1}$ (red points),
$Q_{200,2}$ (black points),
for the collection of functions $[1,f_1,f_2]$, where
$f_1=\log((0.1+1/z)/(1+1/z))$,
$f_2=\log((0.1-1/z)/(1-1/z))$.
}
\label{Fig_log(pL1s)1500(200)_full}
\end{figure}

\begin{figure}[!ht]
\centerline{
\inclgraph{FIG3_2}}
\vskip-6mm
\caption{Numerical distribution of zeros of type~I Hermite--Pad\'e polynomials
$Q_{200,0}$ (blue points)
for the collection of functions $[1,f_1,f_2]$, where
$f_1=\log((0.1+1/z)/(1+1/z))$,
$f_2=\log((0.1-1/z)/(1-1/z))$.
}
\label{Fig_log(pL1s)1500(200)_blue}
\end{figure}

\begin{figure}[!ht]
\centerline{
\inclgraph{FIG3_3}}
\vskip-6mm
\caption{Numerical distribution of zeros of type~I Hermite--Pad\'e polynomials
$P_{200,0}$ (blue points),
$P_{200,1}$ (red points),
$P_{200,2}$ (black points),
for the collection of functions $[1,g_1,g_2]$, where
$g_1=((0.1+1/z)/(1+1/z))^{1/2}$,
$g_2=((0.1-1/z)/(1-1/z))^{1/2}$.
}
\label{Fig_sqrt(pL1s)2000(200)_full}
\end{figure}

\begin{figure}[!ht]
\centerline{
\inclgraph{FIG3_4}}
\vskip-6mm
\caption{Numerical distribution of zeros of type~I Hermite--Pad\'e polynomials
$P_{200,0}$ (blue points)
for the collection of functions $[1,g_1,g_2]$, where
$g_1=((0.1+1/z)/(1+1/z))^{1/2}$,
$g_2=((0.1-1/z)/(1-1/z))^{1/2}$.
}
\label{Fig_sqrt(pL1s)2000(200)_blue}
\end{figure}

%%%endfigures3

\clearpage

%%%figures4

% Figures 13 - 24 (12 total)

\begin{figure}[!ht]
\centerline{
\inclgraph{FIG4_1}}
\vskip-6mm
\caption{Numerical distribution of zeros of type~I Hermite--Pad\'e polynomials
$Q_{200,0}$ (blue points),
$Q_{200,1}$ (red points),
$Q_{200,2}$ (black points),
for the collection of functions $[1,f_1,f_2]$, where
$f_1=\log((0.2-1/z)/(1+1/z))$,
$f_2=\log((0.2+1/z)/(1-1/z))$.
}
\label{Fig_log(p2s)2000(200)_full}
\end{figure}

\begin{figure}[!ht]
\centerline{
\inclgraph{FIG4_2}}
\vskip-6mm
\caption{Numerical distribution of zeros of type~I Hermite--Pad\'e polynomials
$Q_{200,0}$ (blue points)
for the collection of functions $[1,f_1,f_2]$, where
$f_1=\log((0.2-1/z)/(1+1/z))$,
$f_2=\log((0.2+1/z)/(1-1/z))$.
}
\label{Fig_log(p2s)2000(200)_blue}
\end{figure}

\begin{figure}[!ht]
\centerline{
\inclgraph{FIG4_3}}
\vskip-6mm
\caption{Numerical distribution of zeros of type~I Hermite--Pad\'e polynomials
$Q_{200,1}$ (red points)
for the collection of functions $[1,f_1,f_2]$, where
$f_1=\log((0.2-1/z)/(1+1/z))$,
$f_2=\log((0.2+1/z)/(1-1/z))$.
}
\label{Fig_log(p2s)2000(200)_red}
\end{figure}

\begin{figure}[!ht]
\centerline{
\inclgraph{FIG4_4}}
\vskip-6mm
\caption{Numerical distribution of zeros of type~I Hermite--Pad\'e polynomials
$Q_{200,2}$ (black points)
for the collection of functions $[1,f_1,f_2]$, where
$f_1=\log((0.2-1/z)/(1+1/z))$,
$f_2=\log((0.2+1/z)/(1-1/z))$.
}
\label{Fig_log(p2s)2000(200)_bk}
\end{figure}

\begin{figure}[!ht]
\centerline{
\inclgraph{FIG4_5}}
\vskip-6mm
\caption{Numerical distribution of zeros of type~I Hermite--Pad\'e polynomials
$P_{200,0}$ (blue points),
$P_{200,1}$ (red points),
$P_{200,2}$ (black points),
for the collection of functions $[1,g_1,g_2]$, where
$g_1=((0.2-1/z)/(1+1/z))^{1/2}$,
$g_2=((0.2+1/z)/(1-1/z))^{1/2}$.
}
\label{Fig_sqrt(p2s)2000(200)_full}
\end{figure}

\begin{figure}[!ht]
\centerline{
\inclgraph{FIG4_6}}
\vskip-6mm
\caption{Numerical distribution of zeros of type~I Hermite--Pad\'e polynomials
$P_{200,0}$ (blue points)
for the collection of functions $[1,g_1,g_2]$, where
$g_1=((0.2-1/z)/(1+1/z))^{1/2}$,
$g_2=((0.2+1/z)/(1-1/z))^{1/2}$.
}
\label{Fig_sqrt(p2s)2000(200)_blue}
\end{figure}

\begin{figure}[!ht]
\centerline{
\inclgraph{FIG4_7}}
\vskip-6mm
\caption{Numerical distribution of zeros of type~I Hermite--Pad\'e polynomials
$P_{200,1}$ (red points)
for the collection of functions $[1,g_1,g_2]$, where
$g_1=((0.2-1/z)/(1+1/z))^{1/2}$,
$g_2=((0.2+1/z)/(1-1/z))^{1/2}$.
}
\label{Fig_sqrt(p2s)2000(200)_red}
\end{figure}

\begin{figure}[!ht]
\centerline{
\inclgraph{FIG4_8}}
\vskip-6mm
\caption{Numerical distribution of zeros of type~I Hermite--Pad\'e polynomials
$P_{200,2}$ (black points)
for the collection of functions $[1,g_1,g_2]$, where
$g_1=((0.2-1/z)/(1+1/z))^{1/2}$,
$g_2=((0.2+1/z)/(1-1/z))^{1/2}$.
}
\label{Fig_sqrt(p2s)2000(200)_bk}
\end{figure}

\begin{figure}[!ht]
\centerline{
\inclgraph{FIG4_9}}
\vskip-6mm
\caption{Numerical distribution of zeros of type~I Hermite--Pad\'e polynomials
$\UU_{200,0}$ (blue points),
$\UU_{200,1}$ (red points),
$\UU_{200,2}$ (black points),
for the collection of functions $[1,h_1,h_2]$, where
$h_1=((0.2-1/z)/(1+1/z))^{1/3}$,
$h_2=((0.2+1/z)/(1-1/z))^{1/3}$.
}
\label{Fig_cub(p2s)2000(200)_full}
\end{figure}

\begin{figure}[!ht]
\centerline{
\inclgraph{FIG4_10}}
\vskip-6mm
\caption{Numerical distribution of zeros of type~I Hermite--Pad\'e polynomials
$\UU_{200,0}$ (blue points)
for the collection of functions $[1,h_1,h_2]$, where
$h_1=((0.2-1/z)/(1+1/z))^{1/3}$,
$h_2=((0.2+1/z)/(1-1/z))^{1/3}$.
}
\label{Fig_cub(p2s)2000(200)_blue}
\end{figure}

\begin{figure}[!ht]
\centerline{
\inclgraph{FIG4_11}}
\vskip-6mm
\caption{Numerical distribution of zeros of type~I Hermite--Pad\'e polynomials
$\UU_{200,1}$ (red points)
for the collection of functions $[1,h_1,h_2]$, where
$h_1=((0.2-1/z)/(1+1/z))^{1/3}$,
$h_2=((0.2+1/z)/(1-1/z))^{1/3}$.
}
\label{Fig_cub(p2s)2000(200)_red}
\end{figure}

\begin{figure}[!ht]
\centerline{
\inclgraph{FIG4_12}}
\vskip-6mm
\caption{Numerical distribution of zeros of type~I Hermite--Pad\'e polynomials
$\UU_{200,2}$ (black points)
for the collection of functions $[1,h_1,h_2]$, where
$h_1=((0.2-1/z)/(1+1/z))^{1/3}$,
$h_2=((0.2+1/z)/(1-1/z))^{1/3}$.
}
\label{Fig_cub(p2s)2000(200)_bk}
\end{figure}

%%%endfigures4

\clearpage

%%%figures5

% Figures 25 - 36 (12 total)

\begin{figure}[!ht]
\centerline{
\inclgraph{FIG5_1}}
\vskip-6mm
\caption{Numerical distribution of zeros of type~I Hermite--Pad\'e polynomials
$Q_{200,0}$ (blue points),
$Q_{200,1}$ (red points),
$Q_{200,2}$ (black points),
for the collection of functions $[1,f_1,f_2]$, where
$f_1=\log((0.4-1/z)/(1+1/z))$,
$f_2=\log((0.4+1/z)/(1-1/z))$.
}
\label{Fig_log(p4s)2000(200)_full}
\end{figure}

\begin{figure}[!ht]
\centerline{
%\inclgraph{Pic/buslaev_compact}}
\inclgraph{FIG5_2}}
\vskip-6mm
\caption{Numerical distribution of zeros of type~I Hermite--Pad\'e polynomials
$Q_{200,0}$ (blue points)
for the collection of functions $[1,f_1,f_2]$, where
$f_1=\log((0.4-1/z)/(1+1/z))$,
$f_2=\log((0.4+1/z)/(1-1/z))$.
}
\label{Fig_log(p4s)2000(200)_blue}
\end{figure}

\begin{figure}[!ht]
\centerline{
\inclgraph{FIG5_3}}
\vskip-6mm
\caption{Numerical distribution of zeros of type~I Hermite--Pad\'e polynomials
$Q_{200,1}$ (red points)
for the collection of functions $[1,f_1,f_2]$, where
$f_1=\log((0.4-1/z)/(1+1/z))$,
$f_2=\log((0.4+1/z)/(1-1/z))$.
}
\label{Fig_log(p4s)2000(200)_red}
\end{figure}

\begin{figure}[!ht]
\centerline{
\inclgraph{FIG5_4}}
\vskip-6mm
\caption{Numerical distribution of zeros of type~I Hermite--Pad\'e polynomials
$Q_{200,2}$ (black points)
for the collection of functions $[1,f_1,f_2]$, where
$f_1=\log((0.4-1/z)/(1+1/z))$,
$f_2=\log((0.4+1/z)/(1-1/z))$.
}
\label{Fig_log(p4s)2000(200)_bk}
\end{figure}

\begin{figure}[!ht]
\centerline{
\inclgraph{FIG5_5}}
\vskip-6mm
\caption{Numerical distribution of zeros of type~I Hermite--Pad\'e polynomials
$P_{200,0}$ (blue points),
$P_{200,1}$ (red points),
$P_{200,2}$ (black points),
for the collection of functions $[1,g_1,g_2]$, where
$g_1=((0.4-1/z)/(1+1/z))^{1/2}$,
$g_2=((0.4+1/z)/(1-1/z))^{1/2}$.
}
\label{Fig_sqrt(p4s)2000(200)_full}
\end{figure}

\begin{figure}[!ht]
\centerline{
%\inclgraph{Pic/buslaev_compact}}
\inclgraph{FIG5_6}}
\vskip-6mm
\caption{Numerical distribution of zeros of type~I Hermite--Pad\'e polynomials
$P_{200,0}$ (blue points)
for the collection of functions $[1,g_1,g_2]$, where
$g_1=((0.4-1/z)/(1+1/z))^{1/2}$,
$g_2=((0.4+1/z)/(1-1/z))^{1/2}$.
}
\label{Fig_sqrt(p4s)2000(200)_blue}
\end{figure}

\begin{figure}[!ht]
\centerline{
\inclgraph{FIG5_7}}
\vskip-6mm
\caption{Numerical distribution of zeros of type~I Hermite--Pad\'e polynomials
$P_{200,1}$ (red points)
for the collection of functions $[1,g_1,g_2]$, where
$g_1=((0.4-1/z)/(1+1/z))^{1/2}$,
$g_2=((0.4+1/z)/(1-1/z))^{1/2}$.
}
\label{Fig_sqrt(p4s)2000(200)_red}
\end{figure}

\begin{figure}[!ht]
\centerline{
\inclgraph{FIG5_8}}
\vskip-6mm
\caption{Numerical distribution of zeros of type~I Hermite--Pad\'e polynomials
$P_{200,2}$ (black points)
for the collection of functions $[1,g_1,g_2]$, where
$g_1=((0.4-1/z)/(1+1/z))^{1/2}$,
$g_2=((0.4+1/z)/(1-1/z))^{1/2}$.
}
\label{Fig_sqrt(p4s)2000(200)_bk}
\end{figure}

\begin{figure}[!ht]
\centerline{
%\inclgraph{Pic/buslaev_compact}}
\inclgraph{FIG5_9}}
\vskip-6mm
\caption{Numerical distribution of zeros of type~I Hermite--Pad\'e polynomials
$\UU_{200,0}$ (blue points),
$\UU_{200,1}$ (red points),
$\UU_{200,2}$ (black points),
for the collection of functions $[1,h_1,h_2]$, where
$h_1=((0.4-1/z)/(1+1/z))^{1/3}$,
$h_2=((0.4+1/z)/(1-1/z))^{1/3}$.
}
\label{Fig_cub(p4s)2000(200)_full}
\end{figure}

\begin{figure}[!ht]
\centerline{
%\inclgraph{Pic/buslaev_compact}}
\inclgraph{FIG5_10}}
\vskip-6mm
\caption{Numerical distribution of zeros of type~I Hermite--Pad\'e polynomials
$Q_{200,0}$ (blue points)
for the collection of functions $[1,h_1,h_2]$, where
$h_1=((0.4-1/z)/(1+1/z))^{1/3}$,
$h_2=((0.4+1/z)/(1-1/z))^{1/3}$.
}
\label{Fig_cub(p4s)2000(200)_blue}
\end{figure}

\begin{figure}[!ht]
\centerline{
\inclgraph{FIG5_11}}
\vskip-6mm
\caption{Numerical distribution of zeros of type~I Hermite--Pad\'e polynomials
$\UU_{200,1}$ (red points)
for the collection of functions $[1,h_1,h_2]$, where
$h_1=((0.4-1/z)/(1+1/z))^{1/3}$,
$h_2=((0.4+1/z)/(1-1/z))^{1/3}$.
}
\label{Fig_cub(p4s)2000(200)_red}
\end{figure}

\begin{figure}[!ht]
\centerline{
\inclgraph{FIG5_12}}
\vskip-6mm
\caption{Numerical distribution of zeros of type~I Hermite--Pad\'e polynomials
$\UU_{200,2}$ (black points)
for the collection of functions $[1,h_1,h_2]$, where
$h_1=((0.4-1/z)/(1+1/z))^{1/3}$,
$h_2=((0.4+1/z)/(1-1/z))^{1/3}$.
}
\label{Fig_cub(p4s)2000(200)_bk}
\end{figure}

%%%endfigures5

\clearpage

%%%figures6

% Figures 37 - 48 (12 total)

\begin{figure}[!ht]
\centerline{
%\inclgraph{Pic/buslaev_compact}}
\inclgraph{FIG6_1}}
\vskip-6mm
\caption{Numerical distribution of zeros of type~I Hermite--Pad\'e polynomials
$Q_{200,0}$ (blue points),
$Q_{200,1}$ (red points),
$Q_{200,2}$ (black points),
for the collection of functions $[1,f_1,f_2]$, where
$f_1=\log((0.625-1/z)/(1+1/z))$,
$f_2=\log((0.625+1/z)/(1-1/z))$.
}
\label{Fig_log(p625s)2000(200)_full}
\end{figure}

\begin{figure}[!ht]
\centerline{
%\inclgraph{Pic/buslaev_compact}}
\inclgraph{FIG6_2}}
\vskip-6mm
\caption{Numerical distribution of zeros of type~I Hermite--Pad\'e polynomials
$Q_{200,0}$ (blue points)
for the collection of functions $[1,f_1,f_2]$, where
$f_1=\log((0.625-1/z)/(1+1/z))$,
$f_2=\log((0.625+1/z)/(1-1/z))$.
}
\label{Fig_log(p625s)2000(200)_blue}
\end{figure}

\begin{figure}[!ht]
\centerline{
\inclgraph{FIG6_3}}
\vskip-6mm
\caption{Numerical distribution of zeros of type~I Hermite--Pad\'e polynomials
$Q_{200,1}$ (red points)
for the collection of functions $[1,f_1,f_2]$, where
$f_1=\log((0.625-1/z)/(1+1/z))$,
$f_2=\log((0.625+1/z)/(1-1/z))$.
}
\label{Fig_log(p625s)2000(200)_red}
\end{figure}

\begin{figure}[!ht]
\centerline{
\inclgraph{FIG6_4}}
\vskip-6mm
\caption{Numerical distribution of zeros of type~I Hermite--Pad\'e polynomials
$Q_{200,2}$ (black points)
for the collection of functions $[1,f_1,f_2]$, where
$f_1=\log((0.625-1/z)/(1+1/z))$,
$f_2=\log((0.625+1/z)/(1-1/z))$.
}
\label{Fig_log(p625s)2000(200)_bk}
\end{figure}

\begin{figure}[!ht]
\centerline{
%\inclgraph{Pic/buslaev_compact}}
\inclgraph{FIG6_5}}
\vskip-6mm
\caption{Numerical distribution of zeros of type~I Hermite--Pad\'e polynomials
$P_{200,0}$ (blue points),
$P_{200,1}$ (red points),
$P_{200,2}$ (black points),
for the collection of functions $[1,g_1,g_2]$, where
$g_1=((0.625-1/z)/(1+1/z))^{1/2}$,
$g_2=((0.625+1/z)/(1-1/z))^{1/2}$.
}
\label{Fig_sqrt(p625s)2000(200)_full}
\end{figure}

\begin{figure}[!ht]
\centerline{
%\inclgraph{Pic/buslaev_compact}}
\inclgraph{FIG6_6}}
\vskip-6mm
\caption{Numerical distribution of zeros of type~I Hermite--Pad\'e polynomials
$P_{200,0}$ (blue points)
for the collection of functions $[1,g_1,g_2]$, where
$g_1=((0.625-1/z)/(1+1/z))^{1/2}$,
$g_2=((0.625+1/z)/(1-1/z))^{1/2}$.
}
\label{Fig_sqrt(p625s)2000(200)_blue}
\end{figure}

\begin{figure}[!ht]
\centerline{
\inclgraph{FIG6_7}}
\vskip-6mm
\caption{Numerical distribution of zeros of type~I Hermite--Pad\'e polynomials
$P_{200,1}$ (red points)
for the collection of functions $[1,g_1,g_2]$, where
$g_1=((0.625-1/z)/(1+1/z))^{1/2}$,
$g_2=((0.625+1/z)/(1-1/z))^{1/2}$.
}
\label{Fig_sqrt(p625s)2000(200)_red}
\end{figure}

\begin{figure}[!ht]
\centerline{
\inclgraph{FIG6_8}}
\vskip-6mm
\caption{Numerical distribution of zeros of type~I Hermite--Pad\'e polynomials
$P_{200,2}$ (black points)
for the collection of functions $[1,g_1,g_2]$, where
$g_1=((0.625-1/z)/(1+1/z))^{1/2}$,
$g_2=((0.625+1/z)/(1-1/z))^{1/2}$.
}
\label{Fig_sqrt(p625s)2000(200)_bk}
\end{figure}

\begin{figure}[!ht]
\centerline{
%\inclgraph{Pic/buslaev_compact}}
\inclgraph{FIG6_9}}
\vskip-6mm
\caption{Numerical distribution of zeros of type~I Hermite--Pad\'e polynomials
$\UU_{200,0}$ (blue points),
$\UU_{200,1}$ (red points),
$\UU_{200,2}$ (black points),
for the collection of functions $[1,h_1,h_2]$, where
$h_1=((0.625-1/z)/(1+1/z))^{1/3}$,
$h_2=((0.625+1/z)/(1-1/z))^{1/3}$.
}
\label{Fig_cub(p625s)2000(200)_full}
\end{figure}

\begin{figure}[!ht]
\centerline{
%\inclgraph{Pic/buslaev_compact}}
\inclgraph{FIG6_10}}
\vskip-6mm
\caption{Numerical distribution of zeros of type~I Hermite--Pad\'e polynomials
$\UU_{200,0}$ (blue points)
for the collection of functions $[1,h_1,h_2]$, where
$h_1=((0.625-1/z)/(1+1/z))^{1/3}$,
$h_2=((0.625+1/z)/(1-1/z))^{1/3}$.
}
\label{Fig_cub(p625s)2000(200)_blue}
\end{figure}

\begin{figure}[!ht]
\centerline{
\inclgraph{FIG6_11}}
\vskip-6mm
\caption{Numerical distribution of zeros of type~I Hermite--Pad\'e polynomials
$\UU_{200,1}$ (red points)
for the collection of functions $[1,h_1,h_2]$, where
$h_1=((0.625-1/z)/(1+1/z))^{1/3}$,
$h_2=((0.625+1/z)/(1-1/z))^{1/3}$.
}
\label{Fig_cub(p625s)2000(200)_red}
\end{figure}

\begin{figure}[!ht]
\centerline{
\inclgraph{FIG6_12}}
\vskip-6mm
\caption{Numerical distribution of zeros of type~I Hermite--Pad\'e polynomials
$\UU_{200,2}$ (black points)
for the collection of functions $[1,h_1,h_2]$, where
$h_1=((0.625-1/z)/(1+1/z))^{1/3}$,
$h_2=((0.625+1/z)/(1-1/z))^{1/3}$.
}
\label{Fig_cub(p625s)2000(200)_bk}
\end{figure}

%%%endfigures6

\clearpage

%%%figures7

% Figures 49 - 60 (12 total)

\begin{figure}[!ht]
\centerline{
%\inclgraph{Pic/buslaev_compact}}
\inclgraph{FIG7_1}}
\vskip-6mm
\caption{Numerical distribution of zeros of type~I Hermite--Pad\'e polynomials
$Q_{200,0}$ (blue points),
$Q_{200,1}$ (red points),
$Q_{200,2}$ (black points),
for the collection of functions $[1,f_1,f_2]$, where
$f_1=\log((0.73-1/z)/(1+1/z))$,
$f_2=\log((0.73+1/z)/(1-1/z))$.
}
\label{Fig_log(p73s)2000(200)_full}
\end{figure}

\begin{figure}[!ht]
\centerline{
%\inclgraph{Pic/buslaev_compact}}
\inclgraph{FIG7_2}}
\vskip-6mm
\caption{Numerical distribution of zeros of type~I Hermite--Pad\'e polynomials
$Q_{200,0}$ (blue points)
for the collection of functions $[1,f_1,f_2]$, where
$f_1=\log((0.73-1/z)/(1+1/z))$,
$f_2=\log((0.73+1/z)/(1-1/z))$.
}
\label{Fig_log(p73s)2000(200)_blue}
\end{figure}

\begin{figure}[!ht]
\centerline{
\inclgraph{FIG7_3}}
\vskip-6mm
\caption{Numerical distribution of zeros of type~I Hermite--Pad\'e polynomials
$Q_{200,1}$ (red points)
for the collection of functions $[1,f_1,f_2]$, where
$f_1=\log((0.73-1/z)/(1+1/z))$,
$f_2=\log((0.73+1/z)/(1-1/z))$.
}
\label{Fig_log(p73s)2000(200)_red}
\end{figure}

\begin{figure}[!ht]
\centerline{
\inclgraph{FIG7_4}}
\vskip-6mm
\caption{Numerical distribution of zeros of type~I Hermite--Pad\'e polynomials
$Q_{200,2}$ (black points)
for the collection of functions $[1,f_1,f_2]$, where
$f_1=\log((0.73-1/z)/(1+1/z))$,
$f_2=\log((0.73+1/z)/(1-1/z))$.
}
\label{Fig_log(p73s)2000(200)_bk}
\end{figure}

\begin{figure}[!ht]
\centerline{
%\inclgraph{Pic/buslaev_compact}}
\inclgraph{FIG7_5}}
\vskip-6mm
\caption{Numerical distribution of zeros of type~I Hermite--Pad\'e polynomials
$P_{200,0}$ (blue points),
$P_{200,1}$ (red points),
$P_{200,2}$ (black points),
for the collection of functions $[1,g_1,g_2]$, where
$g_1=((0.73-1/z)/(1+1/z))^{1/2}$,
$g_2=((0.73+1/z)/(1-1/z))^{1/2}$.
}
\label{Fig_sqrt(p73s)2000(200)_full}
\end{figure}

\begin{figure}[!ht]
\centerline{
%\inclgraph{Pic/buslaev_compact}}
\inclgraph{FIG7_6}}
\vskip-6mm
\caption{Numerical distribution of zeros of type~I Hermite--Pad\'e polynomials
$P_{200,0}$ (blue points)
for the collection of functions $[1,g_1,g_2]$, where
$g_1=((0.73-1/z)/(1+1/z))^{1/2}$,
$g_2=((0.73+1/z)/(1-1/z))^{1/2}$.
}
\label{Fig_sqrt(p73s)2000(200)_blue}
\end{figure}

\begin{figure}[!ht]
\centerline{
\inclgraph{FIG7_7}}
\vskip-6mm
\caption{Numerical distribution of zeros of type~I Hermite--Pad\'e polynomials
$P_{200,1}$ (red points)
for the collection of functions $[1,g_1,g_2]$, where
$g_1=((0.73-1/z)/(1+1/z))^{1/2}$,
$g_2=((0.73+1/z)/(1-1/z))^{1/2}$.
}
\label{Fig_sqrt(p73s)2000(200)_red}
\end{figure}

\begin{figure}[!ht]
\centerline{
\inclgraph{FIG7_8}}
\vskip-6mm
\caption{Numerical distribution of zeros of type~I Hermite--Pad\'e polynomials
$P_{200,2}$ (black points)
for the collection of functions $[1,g_1,g_2]$, where
$g_1=((0.73-1/z)/(1+1/z))^{1/2}$,
$g_2=((0.73+1/z)/(1-1/z))^{1/2}$.
}
\label{Fig_sqrt(p73s)2000(200)_bk}
\end{figure}

\begin{figure}[!ht]
\centerline{
%\inclgraph{Pic/buslaev_compact}}
\inclgraph{FIG7_9}}
\vskip-6mm
\caption{Numerical distribution of zeros of type~I Hermite--Pad\'e polynomials
$\UU_{300,0}$ (blue points),
$\UU_{300,1}$ (red points),
$\UU_{300,2}$ (black points),
for the collection of functions $[1,h_1,h_2]$, where
$h_1=((0.73-1/z)/(1+1/z))^{1/3}$,
$h_2=((0.73+1/z)/(1-1/z))^{1/3}$.
}
\label{Fig_cub(p73s)1500(300)_full}
\end{figure}

\begin{figure}[!ht]
\centerline{
%\inclgraph{Pic/buslaev_compact}}
\inclgraph{FIG7_10}}
\vskip-6mm
\caption{Numerical distribution of zeros of type~I Hermite--Pad\'e polynomials
$\UU_{300,0}$ (blue points)
for the collection of functions $[1,h_1,h_2]$, where
$h_1=((0.73-1/z)/(1+1/z))^{1/3}$,
$h_2=((0.73+1/z)/(1-1/z))^{1/3}$.
}
\label{Fig_cub(p73s)1500(300)_blue}
\end{figure}

\begin{figure}[!ht]
\centerline{
\inclgraph{FIG7_11}}
\vskip-6mm
\caption{Numerical distribution of zeros of type~I Hermite--Pad\'e polynomials
$\UU_{200,1}$ (red points)
for the collection of functions $[1,h_1,h_2]$, where
$h_1=((0.73-1/z)/(1+1/z))^{1/3}$,
$h_2=((0.73+1/z)/(1-1/z))^{1/3}$.
}
\label{Fig_cub(p73s)1500(300)_red}
\end{figure}

\begin{figure}[!ht]
\centerline{
\inclgraph{FIG7_12}}
\vskip-6mm
\caption{Numerical distribution of zeros of type~I Hermite--Pad\'e polynomials
$\UU_{200,2}$ (black points)
for the collection of functions $[1,h_1,h_2]$, where
$h_1=((0.73-1/z)/(1+1/z))^{1/3}$,
$h_2=((0.73+1/z)/(1-1/z))^{1/3}$.
}
\label{Fig_cub(p73s)1500(300)_bk}
\end{figure}

%%%endfigures7

\clearpage

%%%figures8

% Figures 61 - 72 (12 total)

\begin{figure}[!ht]
\centerline{
%\inclgraph{Pic/buslaev_compact}}
\inclgraph{FIG8_1}}
\vskip-6mm
\caption{Numerical distribution of zeros of type~I Hermite--Pad\'e polynomials
$Q_{200,0}$ (blue points),
$Q_{200,1}$ (red points),
$Q_{200,2}$ (black points),
for the collection of functions $[1,f_1,f_2]$, where
$f_1=\log((0.8-1/z)/(1+1/z))$,
$f_2=\log((0.8+1/z)/(1-1/z))$.
}
\label{Fig_log(p8s)2000(200)_full}
\end{figure}

\begin{figure}[!ht]
\centerline{
%\inclgraph{Pic/buslaev_compact}}
\inclgraph{FIG8_2}}
\vskip-6mm
\caption{Numerical distribution of zeros of type~I Hermite--Pad\'e polynomials
$Q_{200,0}$ (blue points)
for the collection of functions $[1,f_1,f_2]$, where
$f_1=\log((0.8-1/z)/(1+1/z))$,
$f_2=\log((0.8+1/z)/(1-1/z))$.
}
\label{Fig_log(p8s)2000(200)_blue}
\end{figure}

\begin{figure}[!ht]
\centerline{
\inclgraph{FIG8_3}}
\vskip-6mm
\caption{Numerical distribution of zeros of type~I Hermite--Pad\'e polynomials
$Q_{200,1}$ (red points)
for the collection of functions $[1,f_1,f_2]$, where
$f_1=\log((0.8-1/z)/(1+1/z))$,
$f_2=\log((0.8+1/z)/(1-1/z))$.
}
\label{Fig_log(p8s)2000(200)_red}
\end{figure}

\begin{figure}[!ht]
\centerline{
\inclgraph{FIG8_4}}
\vskip-6mm
\caption{Numerical distribution of zeros of type~I Hermite--Pad\'e polynomials
$Q_{200,2}$ (black points)
for the collection of functions $[1,f_1,f_2]$, where
$f_1=\log((0.8-1/z)/(1+1/z))$,
$f_2=\log((0.8+1/z)/(1-1/z))$.
}
\label{Fig_log(p8s)2000(200)_bk}
\end{figure}

\begin{figure}[!ht]
\centerline{
%\inclgraph{Pic/buslaev_compact}}
\inclgraph{FIG8_5}}
\vskip-6mm
\caption{Numerical distribution of zeros of type~I Hermite--Pad\'e polynomials
$P_{200,0}$ (blue points),
$P_{200,1}$ (red points),
$P_{200,2}$ (black points),
for the collection of functions $[1,g_1,g_2]$, where
$g_1=((0.8-1/z)/(1+1/z))^{1/2}$,
$g_2=((0.8+1/z)/(1-1/z))^{1/2}$.
}
\label{Fig_sqrt(p8s)2000(200)_full}
\end{figure}

\begin{figure}[!ht]
\centerline{
%\inclgraph{Pic/buslaev_compact}}
\inclgraph{FIG8_6}}
\vskip-6mm
\caption{Numerical distribution of zeros of type~I Hermite--Pad\'e polynomials
$P_{200,0}$ (blue points)
for the collection of functions $[1,g_1,g_2]$, where
$g_1=((0.8-1/z)/(1+1/z))^{1/2}$,
$g_2=((0.8+1/z)/(1-1/z))^{1/2}$.
}
\label{Fig_sqrt(p8s)2000(200)_blue}
\end{figure}

\begin{figure}[!ht]
\centerline{
\inclgraph{FIG8_7}}
\vskip-6mm
\caption{Numerical distribution of zeros of type~I Hermite--Pad\'e polynomials
$P_{200,1}$ (red points)
for the collection of functions $[1,g_1,g_2]$, where
$g_1=((0.8-1/z)/(1+1/z))^{1/2}$,
$g_2=((0.8+1/z)/(1-1/z))^{1/2}$.
}
\label{Fig_sqrt(p8s)2000(200)_red}
\end{figure}

\begin{figure}[!ht]
\centerline{
\inclgraph{FIG8_8}}
\vskip-6mm
\caption{Numerical distribution of zeros of type~I Hermite--Pad\'e polynomials
$P_{200,2}$ (black points)
for the collection of functions $[1,g_1,g_2]$, where
$g_1=((0.8-1/z)/(1+1/z))^{1/2}$,
$g_2=((0.8+1/z)/(1-1/z))^{1/2}$.
}
\label{Fig_sqrt(p8s)2000(200)_bk}
\end{figure}

\begin{figure}[!ht]
\centerline{
%\inclgraph{Pic/buslaev_compact}}
\inclgraph{FIG8_9}}
\vskip-6mm
\caption{Numerical distribution of zeros of type~I Hermite--Pad\'e polynomials
$\UU_{200,0}$ (blue points),
$\UU_{200,1}$ (red points),
$\UU_{200,2}$ (black points),
for the collection of functions $[1,h_1,h_2]$, where
$h_1=((0.85-1/z)/(1+1/z))^{1/3}$,
$h_2=((0.85+1/z)/(1-1/z))^{1/3}$.
}
\label{Fig_cub(p8s)2000(200)_full}
\end{figure}

\begin{figure}[!ht]
\centerline{
%\inclgraph{Pic/buslaev_compact}}
\inclgraph{FIG8_10}}
\vskip-6mm
\caption{Numerical distribution of zeros of type~I Hermite--Pad\'e polynomials
$\UU_{200,0}$ (blue points)
for the collection of functions $[1,h_1,h_2]$, where
$h_1=((0.85-1/z)/(1+1/z))^{1/3}$,
$h_2=((0.85+1/z)/(1-1/z))^{1/3}$.
}
\label{Fig_cub(p8s)2000(200)_blue}
\end{figure}

\begin{figure}[!ht]
\centerline{
\inclgraph{FIG8_11}}
\vskip-6mm
\caption{Numerical distribution of zeros of type~I Hermite--Pad\'e polynomials
$\UU_{200,1}$ (red points)
for the collection of functions $[1,h_1,h_2]$, where
$h_1=((0.85-1/z)/(1+1/z))^{1/3}$,
$h_2=((0.85+1/z)/(1-1/z))^{1/3}$.
}
\label{Fig_cub(p8s)2000(200)_red}
\end{figure}

\begin{figure}[!ht]
\centerline{
\inclgraph{FIG8_12}}
\vskip-6mm
\caption{Numerical distribution of zeros of type~I Hermite--Pad\'e polynomials
$\UU_{200,2}$ (black points)
for the collection of functions $[1,h_1,h_2]$, where
$h_1=((0.85-1/z)/(1+1/z))^{1/3}$,
$h_2=((0.85+1/z)/(1-1/z))^{1/3}$.
}
\label{Fig_cub(p8s)2000(200)_bk}
\end{figure}

%%%endfigures8

\end{document}